\documentclass[12pt,twoside,a4paper,english]{smfart}

\usepackage[utf8]{inputenc}
\usepackage[francais,english]{babel}
\usepackage{lmodern}
\usepackage[T1]{fontenc}

\usepackage{amsfonts,amssymb,amsmath}
\usepackage{amsthm}
\usepackage{mathrsfs}
\usepackage{geometry}
\usepackage[all]{xy}

\usepackage[plainpages=true,naturalnames=true,hyperfootnotes=false,colorlinks=false,hypertexnames=false,linkbordercolor={0 1 0},citebordercolor={1 0 0}]{hyperref}

\newcommand{\A}{\mathbf{A}}

\newcommand{\D}{\mathrm{D}}
\newcommand{\R}{\mathbf{R}}
\newcommand{\Z}{\mathbf{Z}}

\newcommand{\F}{\mathcal{F}}

\newcommand{\spec}{\mathrm{Spec}\, }

\newcommand{\proj}{\mathrm{Proj}\, }

\renewcommand{\phi}{\varphi}
\renewcommand{\epsilon}{\varepsilon}
\renewcommand{\H}{\mathrm{H}}

\newcommand{\pt}{\mathrm{pt}}
\newcommand{\id}{\mathrm{id}}
\newcommand{\res}{\mathrm{r{e}s}}
\newcommand{\syl}{\mathrm{Syl}}

\newcommand{\susp}{\mathbf{\Sigma}}

\newcommand{\bez}{\mathrm{B\acute{e}z}}
\newcommand{\SL}{\mathbf{SL}}

\newcommand{\pun}{\mathbf{P}^1}
\newcommand{\aun}{\mathbf{A}^1}
\newcommand{\pd}{\mathbf{P}^d}
\newcommand{\gm}{\mathbf{G_{m}}}
\newcommand{\ga}{\mathbf{G_{a}}}

\newcommand{\fd}{\longrightarrow}
\newcommand{\xfd}{\xrightarrow}
\newcommand{\dfd}{\rightrightarrows}

\newcommand{\rond}{\circ}
\newcommand{\tens}{\otimes}
\newcommand{\plus}{\oplus\naif}

\newcommand{\isom}{\simeq}

\newcommand{\udisj}{\coprod}

\newcommand{\nn}{$n\hspace{-.07cm}\times\hspace{-.07cm}n$}

\renewcommand{\geq}{\geqslant}
\renewcommand{\leq}{\leqslant}

\newcommand{\mc}{\mathcal}

\newcommand*{\croix}[1]{{\underset{#1} {\times}}}

\newcommand*{\sur}[2]{
\raisebox{0,2ex}{{$ #1$}}
\raisebox{-1ex}
{\Large \slash}_{\hspace{-1ex}\raisebox{0,75ex}
{{ $#2$}}} }

\newcommand{\sura}[2]{
\raisebox{0,2ex}{{$#1$}}
\hspace{-1ex}
\raisebox{-1ex}
{\Large \slash}_{\raisebox{0,75ex}
{{$#2$}}} }

\newcommand*{\surkk}[2]{
\raisebox{0,2ex}{\ensuremath{\scriptstyle #1}}
\hspace{-1ex}
\raisebox{-1ex}
{{\slash}}_{\raisebox{0,35ex}
{\ensuremath{ \scriptstyle #2}}} }

\newcommand{\kk}{ \surkk{k^{\times}}{k^{\times 2}}}
\newcommand{\kkn}{ \surkk{k^{\times}}{k^{\times 2n}}}

\newcommand{\app}[5]{ 
  $$
   \begin{array}{lccl}
       #1 :& #2 & \fd & #3\\
      & #4 &\mapsto & #5
    \end{array} 
$$
}

\newcommand{\appsn}[4]{ 
$$
   \begin{array}{ccl}
        #1 & \fd & #2\\
       #3 &\mapsto & #4
    \end{array}
$$
}

\newcommand{\appsnbis}[4]{ 
$$
   \begin{array}{rcl}
        #1 & \fd & #2\\
       #3 &\mapsto & #4
    \end{array}
$$
}

\newcommand{\cache}[1]{}

\newcommand{\inj}{\hookrightarrow}

\renewcommand{\tilde}{\widetilde}

\newcommand{\hocolim}[1]{\underset{#1}{\mathrm{hocolim}}\ }

\newcommand{\cf}{\emph{c.f. }}

\newcommand{\ens}{\mathscr{S}\hspace{-.4ex}et}

\newcommand{\kalg}{\mathrm{Alg}_k}

\newcommand{\ie}{\emph{i.e.}}

\newcommand{\attention}{\emph{Warning: }}

\newcommand{\pointloin}{\qquad.}

\newcommand{\bp}[1]{\big{(}#1 \big{)}}
\newcommand{\Bp}[1]{\Big{(}#1 \Big{)}}

\newcommand{\sph}{\mathrm{S}}

\newcommand{\et}{\quad \text{ et } \quad}
\newcommand{\tr}{\mathrm{tr}} 

\newcommand{\hp}{\stackrel{\mathrm{p}}{\sim}}
\newcommand{\hl}{\stackrel{\mathrm{u}}{\sim}}
\newcommand{\mot}{{\aun}}
\newcommand{\pp}{\left[\pun, \pun \right]}
\newcommand{\ppl}{\left\lbrace \pun, \pun \right\rbrace \naif}
\newcommand{\ppd}{\left[\pun, \pd \right]}

\newcommand{\MW}{\mathrm{MW}}
\newcommand{\GW}{\mathrm{GW}}
\newcommand{\libre}{\mathscr{U}}
\newcommand{\transp}[1]{{\vphantom{#1}}^{\mathit  t}{#1}}
\newcommand{\hank}{\mathrm{Hank}}

\newcommand{\dfxy}[2]{\ar@<-.5ex>[r]_-{#1}\ar@<.5ex>[r]^-{#2}}

\newcommand{\esp}{\mathscr{S}}

\newcommand{\ho}{\mathrm{Ho}}

\newcommand{\naif}{^{\mathrm{N}}}
\newcommand{\ppn}{\pp\naif}
\newcommand{\ppdn}{\ppd\naif}

\newcommand{\pN}{\pi_0\naif}

\newcommand{\tMW}{\MW^{\mathrm{s}}}

\newcommand{\vtruc}{\, {\textstyle \udisj\limits_{0 \sim \infty}}\,}
\newcommand{\vtrucp}[2]{\,{\textstyle\udisj\limits_{#1 \sim #2}}\,}
 
\entrymodifiers={!!<0pt,0.7ex>+}

\newtheorem{theoreme}{Theorem}[section]

\newtheorem{proposition}[theoreme]{Proposition}
\newtheorem{prop-def}[theoreme]{Proposition-definition}
\newtheorem{corollaire}[theoreme]{Corollary}
\newtheorem{lemme}[theoreme]{Lemma}

\theoremstyle{definition}

\newtheorem{definition}[theoreme]{Definition}
\newtheorem{exemple}[theoreme]{Example}
\newtheorem{exemples}[theoreme]{Examples}
\newtheorem{remarque}[theoreme]{Remark}
\newtheorem{remarques}[theoreme]{Remarks}

\renewcommand{\theequation}{\thetheoreme}

\newcommand{\comm}[1]{{\bf {\sc {\Large (#1)}} }} 
\renewcommand{\comm}[1]{} 
\entrymodifiers={!!<0pt,0.7ex>+}
\title{Algebraic homotopy classes of rational functions}
\alttitle{Classes d'homotopie algébrique de fractions rationnelles}
\author{Christophe Cazanave}
\address{Hausdorff center for mathematics\\
Endenicher Allee, 60 \\
Bonn, Germany}
\email{cazanave@math.uni-bonn.de}
\date{\today}

\begin{document}
\frontmatter

\begin{abstract}
 We compute 
the set ${\left[\mathbf{P}^1, \mathbf{P}^1 \right]}^{\mathrm{N}}$ of \emph{naive} homotopy classes 
of scheme endomorphisms 
of the projective line $\mathbf{P}^1$ over the spectrum of a field.
 Our result compares well with Morel's computation in \cite{morel} of the 
\emph{group} ${\left[\mathbf{P}^1, \mathbf{P}^1 \right]}^{\mathbf{A}^1}$ of ${\mathbf{A}^1}$-homotopy classes of endomorphisms of $\mathbf{P}^1$: the set ${\left[\mathbf{P}^1, \mathbf{P}^1\right]}^{\mathrm{N}}$ admits an \emph{a priori} monoid structure such that the canonical map
${\left[\mathbf{P}^1, \mathbf{P}^1 \right]}^{\mathrm{N}} \rightarrow {\left[\mathbf{P}^1, \mathbf{P}^1 \right]}^{\mathbf{A}^1}$
 is a group completion.
\end{abstract}

\selectlanguage{francais}

\begin{altabstract}
Nous déterminons l'ensemble ${\left[\mathbf{P}^1, \mathbf{P}^1 \right]}^{\mathrm{N}}$ 
des classes d'homotopie naïve d'endomorphismes de 
schémas de la droite projective $\mathbf{P}^1$ au-dessus du spectre d'un corps.
Notre résultat se compare bien avec le calcul, effectué par Morel dans
\cite{morel}, du \emph{groupe} ${\left[\mathbf{P}^1, \mathbf{P}^1 \right]}^{\mathbf{A}^1}$ des classes d'homotopie ${\mathbf{A}^1}$
d'endomorphismes de $\mathbf{P}^1$:  l'ensemble ${\left[\mathbf{P}^1, \mathbf{P}^1 \right]}^{\mathrm{N}}$ admet \emph{a priori} une structure
de monoïde  telle que l'application canonique 
${\left[\mathbf{P}^1, \mathbf{P}^1 \right]}^{\mathrm{N}} \rightarrow {\left[\mathbf{P}^1, \mathbf{P}^1 \right]}^{\mathbf{A}^1}$ soit une compl\'etion en groupe.
\end{altabstract}

\selectlanguage{english}

\thanks{This research was partially supported by the project ANR blanc BLAN08-2\_338236, HGRT}

\maketitle

\mainmatter
\section{Introduction}

The work of Fabien Morel and Vladimir Voevodsky on $\aun$-homotopy theory 
\cite{morel_asterisque,MV} provides a convenient 
framework to do algebraic topology 
in the setting of algebraic geometry. 
More precisely, for any fixed field $k$, Morel and Voevodsky defined 
a convenient category of \emph{spaces}, say $\esp$, containing the
 category of smooth algebraic $k$-varieties 
as a full subcategory, which they endowed with a
 suitable model structure,
in the sense of Quillen's homotopical algebra \cite{quillen}. 
Thus, given two  spaces $X$ and $Y$ in $\esp$ (resp. two pointed spaces),
the set $\lbrace X,Y\rbrace^\mot$ (resp. the set $[X,Y]^\mot$)
 of $\aun$-homotopy classes of \emph{unpointed} morphisms (resp. of \emph{pointed} morphisms)
 from $X$ to $Y$ is well defined and has all the properties 
an algebraic topologist can expect. 
However,  for concrete $X$ and $Y$, these sets are in general very hard to compute.\\

At the starting point of $\aun$-homotopy theory is the notion of
\emph{naive homotopy\footnote{In \cite{MV}, the authors use the terminology  ``elementary homotopy''.}} between two morphisms in $\esp$. Its definition, first introduced by Karoubi and Villamayor \cite{KV}, 
mimics that of the usual notion of homotopy between two topological maps,
replacing the unit interval $[0,1]$ by its algebraic analogue, the affine line
$\aun$.

\begin{definition}
\label{homo_naive}
 Let $X$ and $Y$ be two spaces in $\esp$. 
A naive homotopy is a morphism
$$ F: X \times \aun \fd Y \pointloin$$ 
The restriction $\sigma(F):=F_{\vert X \times \{0\}}$ is the source of 
the homotopy and $\tau(F):=F_{\vert X \times \{1\}}$ is its target.\\
When $X$ and $Y$ have base points, say $x_0$ and $y_0$, we say that $F$ is pointed if its restriction to 
 $\{x_0\}\times \aun$ is constant equal to $y_0$.
\end{definition}

With this notion, one defines the set $\lbrace X,Y \rbrace \naif$  (resp. the set $[X,Y]\naif$)
 of unpointed (resp. of pointed) \emph{naive} 
homotopy classes of morphisms from $X$ to $Y$ as the quotient of the set
 of unpointed (resp. of pointed) morphisms by the equivalence 
relation generated
 by unpointed (resp. by pointed) naive homotopies. 
These sets are sometimes easier to compute than their $\aun$ analogues, 
but they are not very well behaved. There is a canonical map
$$ [X,Y]\naif \fd [X,Y]^\mot$$
which in general is far from being a bijection.
This article studies a particular example where this
map has an interesting behaviour.\\

Let $k$ be a base field. In this article, we focus on the set
 of pointed homotopy classes of endomorphisms of the projective line 
$\pun$ over $\spec k$. Here, our convention is to take the base point in $\pun$ at $\infty:=[1:0]$. 
The main result of this article is the following (see theorem~\ref{th_fr}, corollary~\ref{coro_th_effectif} and theorem~\ref{conjecture_comparaison} for
 more precise and explicit statements).

\begin{theoreme}
\label{th_principal}
The set $\ppn$ of pointed naive homotopy classes of endomorphisms of the projective line admits an \emph{a priori} monoid structure, whose law is denoted by $\plus$. The canonical map
from the monoid $(\ppn,\plus)$ 
to the group $(\pp^\mot, \oplus^\mot)$ of $\aun$-homotopy classes of endomorphisms of $\pun$ is a group completion. 
\end{theoreme}

\subsection*{Overview}
We summarize now briefly the organization of the paper.
\begin{itemize}
\item Section~\ref{sec_fr} reviews the classical
 correspondence between scheme endomorphisms of the projective  line over $\spec k$
and rational functions with coefficients in $k$. 
We give also an analogous description
of pointed naive homotopies of endomorphisms of $\pun$
in terms of rational functions
 with coefficients in the polynomial ring $k[T]$.  

\item Section~\ref{sec_coeur} is the core of the article.
\S\ref{loi_plus} describes a natural monoid structure 
on the scheme of pointed rational functions which induces the monoid structure on $\ppn$ mentioned in theorem~\ref{th_principal}. In \S\ref{subsec_bezout}, we 
review a classical construction, due to Bézout, which associates a non-degenerate symmetric bilinear form to any rational function. 
This leads to the definition of a homotopy invariant, which is at the center of our study.
In \S\ref{subsec_enonce},  the main theorem is stated: 
the previous homotopy invariant
distinguishes all the homotopy classes of rational functions. The proof is given  
 in \S\ref{subsec_dem}.
 Finally, \S\ref{subsec_comparaison} compares our result  
to the actual $\aun$-homotopy classes, as computed by Morel in \cite{morel}.
\item In section~\ref{sec_extensions}, we discuss 
some natural extensions of our main result.
We first give an explicit description of the
 \emph{unpointed} naive homotopy classes of endomorphisms of $\pun$ in \S\ref{subsec_libre}.
Next, in \S\ref{subsec_composition}, we study the composition of endomorphisms
of $\pun$ in terms of our description of $\ppn$.
Finally, for every integer
$d\geq 2$, we compute the set $\ppdn$ of pointed naive
homotopy classes of morphisms from $\pun$ to $\pd$. Not surprisingly, the case $d\geq 2$ is much easier
than the case $d=1$. Our result still compares well with
Morel's computation of the group $\ppd^\mot$ of $\aun$-homotopy classes of
morphisms from $\pun$ to $\pd$. 
\end{itemize}
Three appendices conclude the article.
\begin{itemize}
\item  In appendix~\ref{sec_hermite}, we give 
an elementary computation, due to Ojanguren \cite{ojan}, of 
the set of naive homotopy classes of non-degenerate
 symmetric matrices. It is based on an elegant use of Hermite
 inequality for symmetric bilinear forms over the ring $k[T]$.
\item Appendix~\ref{appendice_comparaison} proves in detail 
that the canonical map $\ppn \fd \pp^\mot$ is compatible with
 monoid structures. 
\item  Appendix~\ref{sec_hankel} reviews the classical correspondence
 between  rational functions and non-degenerate Hankel matrices.
\end{itemize}

\subsection*{Acknowledgements}
The material presented here constitutes the first part 
of the author's Ph-D thesis \cite{these}. The main result was 
first announced in the note \cite{crass} when the base field is 
of characteristic not 2. 
  I am very much indebted to Jean Lannes for his precious
 and constant help.

\section{Rational functions and naive homotopies}
\label{sec_fr}

This section reviews the classical correspondence 
between pointed endomorphisms of the projective line $\pun$
(with base point $\infty$)
 over the spectrum of the field $k$ and pointed
rational functions with coefficients in $k$. We also give a concrete description of pointed naive homotopies of  endomorphisms of $\pun$ 
in terms of pointed rational functions with coefficients in the polynomial ring $k[T]$.

\subsection{Pointed endomorphisms of the projective line}
We first introduce some notation.

\begin{definition}
 For every positive integer $n$, the scheme $\F_n$ of
pointed degree $n$ rational functions is the open subscheme of the affine space
$\A^{2n}=\spec(k[a_0,\dots, a_{n-1},b_0,\dots, b_{n-1}])$
complementary to the hypersurface of equation\footnote{The notation $\res_{n,n}(A,B)$ stands for the 
resultant of the polynomials $A$ and $B$ (considered as polynomials of degree less or equal to $n$). Our conventions, in particular for the sign in the Bézout formula~(\ref{bez_formula}) below, are those of
Bourbaki, \cite{Bourbaki}, \S6,~$\mathrm{n}^{\rond}$6,~IV.}
$$\res_{n,n}(X^n+a_{n-1}X^{n-1} + \dots + a_0 ,b_{n-1}X^{n-1} + \dots + b_0)=0 \pointloin$$
\end{definition}

\begin{remarque}
\label{remarque_res_et_Bezout}
Let $R$ be a ring and $n$ a non-negative integer. 
By the very definition, an $R$-point of $\F_n$
is a pair $(A,B)$ of polynomials of $R[X]$, where
\begin{itemize}
\item $A$ is monic 
of degree $n$,
\item  $B$ is of degree strictly less than $n$,
\item the scalar $\res_{n,n}(A,B)$ is invertible in $R$.
\end{itemize}
Such a point is denoted by $\frac{A}{B}$ and is called a pointed degree $n$ rational function with coefficients in $R$.\\
In the sequel, it is useful to remark that 
the above condition $\res_{n,n}(A,B) \in R^\times$
 is equivalent to the existence of a 
(necessarily unique) Bézout relation
$$ AU+BV=1$$
with $U$ and $V$ polynomials in $R[X]$ such that
$\deg U \leq n-2$ and $\deg V \leq n-1$.
\end{remarque}

The precise correspondence between endomorphisms of $\pun$ and  rational functions is  now summarized in the following proposition.

\begin{proposition}
\label{prop_description}
The datum of a pointed scheme endomorphism of 
$\pun$ over $\spec k$, say $f$, is equivalent to the datum
of an non-negative integer $n$ and of an 
element $\frac{A}{B}\in \F_n(k)$.\\
 The integer $n$ is called the \emph{degree} of $f$ and is denoted $\deg(f)$; 
the scalar $\res_{n,n}(A,B)$ is called the \emph{resultant} of $f$ 
and is denoted $\res(f)$.
\end{proposition}

\subsection{Naive homotopies}

Recall from definition~\ref{homo_naive} that a pointed naive homotopy of endomorphisms of $\pun$, is a scheme morphism
$$ F:\pun \times \aun \fd \pun $$ 
satisfying the appropriate base point condition. The following
 slight generalization of
proposition~\ref{prop_description} gives a description
of pointed naive homotopies in terms of rational functions with coefficients in $k[T]$.

\begin{proposition}
\label{prop_decodage_homotopie}
 The datum of a pointed naive homotopy 
$F:\pun \times \aun \fd \pun$ is equivalent to
the datum of a non-negative integer $n$ and of an element
in $\F_n(k[T])$.
 The source $\sigma(F)$ and the target $\tau(F)$ of $F$ are obtained 
by evaluating the indeterminate $T$ at $0$ and $1$ respectively.   
\end{proposition}

\begin{exemple} \label{exemple_homotopies}
Let $n$ be a positive integer and $b_0$ be a unit in $k^\times$. 
\begin{enumerate}
 \item \label{exemple_homotopies_polynome}Let $A=X^n + a_{n-1}X^{n-1} + \dots +a_0$ be a monic degree $n$ polynomial of $k[X]$. \\
The element  
$$\frac{X^n + T a_{n-1}X^{n-1} + \dots + Ta_0}{b_0} \in \F_n(k[T])$$
gives a pointed naive homotopy between
$\frac{A}{b_0}$ and $\frac{X^n}{b_0}$.
 In other terms, any polynomial is homotopic to its leading term.
\item \label{exemple_homotopies_inverse_polynome}Let $B=b_{n-1}X^{n-1} + \dots +b_0$ be a polynomial of degree $\leq n-1$ such that $B(0)=b_0$.
Then $\frac{X^n}{B}$ is a $k$-point of 
$\F_n$ and the element
$$ \frac{X^n}{T b_{n-1}X^{n-1} + \dots +T b_1 X+ b_0} \in \F_n(k[T])$$ 
gives a pointed naive homotopy between 
 $\frac{X^n}{B}$ and $\frac{X^n}{b_0}$.
\end{enumerate}
\end{exemple}

The simplicity of the above examples of homotopies 
is somewhat atypical.
  In general, given a random rational function,
it is not \emph{a priori} easy to find non-trivial homotopies.
In \S\ref{loi_plus}, we  give a way to
 produce some such families (in fact just enough
to be able to determine $\ppn$), see remark~\ref{rem_loi_plus_sur_pi0}\,(\ref{homo_non_triviale}).

\begin{definition} 
\label{def_hp}
Let $f$ and $g$ be two pointed rational functions over $k$.
We say that $f$ and $g$ are in the same pointed naive homotopy class, and we write $f\hp g$,
if there exists a finite sequence of pointed homotopies,  
say $(F_i)$ with $0 \leq i\leq N$, such that  
\begin{itemize}
\item $\sigma(F_0)=f$ and $\tau(F_N)=g$;
\item for every $ 0\leq i \leq N-1$, we have $\tau(F_i)=\sigma(F_{i+1})$. 
\end{itemize}
The set of pointed naive homotopy classes $\ppn$ 
identifies by definition with the quotient of the set 
$\udisj\limits_{n \geq 0} \F_n(k)$ of pointed rational functions with coefficients in $k$ by the equivalence relation $\hp$.
\end{definition}
 
Note that proposition~\ref{prop_decodage_homotopie} implies
that  two pointed rational functions which are in the same pointed naive homotopy class 
  have same degree and same resultant.
In particular, the set $\ppn$ splits
as the disjoint union of its components of a given degree
$$ \ppn= \udisj_{n\geq0} \ppn_n \qquad.$$

\begin{remarque} 
\label{rem_invariants}
Pointed naive homotopies of rational functions are algebraic paths in the scheme
 of pointed rational functions.
It is convenient to reformulate the 
preceding discussion in terms of the ``naive connected 
components'' of this scheme.\\
Let $\mc{G}:\kalg \fd \ens$ be a functor from the category of $k$-algebras to that of sets. Recall that the naive connected components of $\mc{G}$
 is the new functor $\pN\mc{G}:\kalg \fd \ens$ which assigns
to any $k$-algebra $R$ the coequalizer of the double-arrow
$$ \mc{G}(R[T]) \dfd \mc{G}(R)$$
given by evaluation at $T=0$ and $T=1$.
Moreover, any natural transformation $\mathcal{T}:\mc{G} \fd \mc{H}$ between two 
such functors  induces a natural transformation 
$\pN\mc{T}:\pN\mc{G} \fd \pN\mathcal{H}$.
\\
For every non-negative integer $n$, 
proposition~\ref{prop_decodage_homotopie}
gives a bijection
$$ \ppn_n \isom (\pN\F_n)(k) \pointloin$$
Any scheme morphism $\F_n \fd X$ produces a homotopy invariant
$(\pN\F_n)(k) \fd (\pN X)(k)$.
\end{remarque}

\section{Homotopy classes of rational functions}
\label{sec_coeur}
This section is the core of the article. 
We first endow the set $\ppn$ with an \emph{a priori} monoid
structure in \S\ref{subsec_loi_plus}. We then review 
in \S\ref{subsec_bezout} a classical
construction due to Bézout which associates to any rational 
function a non-degenerate symmetric matrix called the 
\emph{Bézout form}. This leads us to the definition of a homotopy invariant.
The main result, theorem~\ref{th_fr}, is stated in \S\ref{subsec_enonce}
and proved in \S\ref{subsec_dem}. It shows
 that the homotopy invariant
associated to the Bézout form distinguishes all the pointed naive 
 homotopy classes of rational functions.
\S\ref{subsec_comparaison} finally concludes the section by comparing our naive result to the 
actual $\aun$-computation, due to Morel. The result is that the set $\pp^\mot$ has for formal reasons a group structure
and it turns out that this group is isomorphic to the group completion of the monoid $(\ppn, \plus)$.

\subsection{Additions of rational functions}
\label{loi_plus} \label{subsec_loi_plus}

An important feature in the statement 
of our results, is the existence of a
 monoid structure
exhibited \emph{a priori} on the set of naive homotopy classes $\ppn$.
 In fact, there is even a 
graded monoid structure on the disjoint union scheme 
$$\F:=\udisj\limits_{n\geq 0} \F_n \pointloin$$
The datum of such a structure is equivalent to 
the datum of a family of morphisms
$\F_{n_1} \times \F_{n_2} \fd \F_{n_1+n_2}$
indexed by pairs $(n_1,n_2)$ of non-negative integers and subject to 
an associativity condition.\\

Let $(n_1,n_2)$ be a pair of non-negative integers. We
now describe the above structural morphism on the level of functors of points, that
is as a natural transformation of functors from the category of $k$-algebras to the category of sets:
$$ \F_{n_1}(-) \times \F_{n_2}(-) \fd \F_{n_1+n_2}(-) \pointloin$$
Let $R$ be a $k$-algebra.  Two rational functions
$\frac{A_i}{B_i}\in \F_{n_i}(R)$, for $i=1,2$, 
uniquely define two pairs $(U_i,V_i)$ of polynomials 
of $R[X]$ with 
  $\deg U_i \leq n_i-2$ and $\deg V_i \leq n_i-1$
and satisfying Bézout identities
$A_i U_i+B_i V_i=1$ (see remark~\ref{remarque_res_et_Bezout}).
 Define polynomials $A_3, B_3, U_3$ and $V_3$ by setting\footnote{The dot in the right-hand term stands for the usual matrix multiplication.}
  $$
\begin{bmatrix}A_3 & -V_3\\ B_3 & U_3\end{bmatrix} :=
\begin{bmatrix}A_1 & -V_1\\ B_1 & U_1\end{bmatrix} \cdot
\begin{bmatrix}A_2 & -V_2\\ B_2 &U_2\end{bmatrix} \qquad.
$$
The matrices $\begin{bmatrix}A_1 & -V_1\\ B_1 & U_1\end{bmatrix}$ and $\begin{bmatrix}A_2 & -V_2\\ B_2 & U_2\end{bmatrix}$ belong to $\SL_2(k[T])$, thus this is also the case for 
$\begin{bmatrix}A_3 & -V_3\\ B_3 & U_3\end{bmatrix}$. This means that
we have also a Bézout relation for the polynomials $A_3$ and $B_3$.
Moreover, observe that 
$A_3=A_1A_2 - V_1B_2$ is monic of degree $n_1+n_2$ and that $B_3=B_1A_2 + U_1 B_2$ is of degree strictly less than $n_1+n_2$. So  $\frac{A_3}{B_3}$ is an $R$-point of $\F_{n_1+n_2}$.
We write 
$$\frac{A_1}{B_1} \plus \frac{A_2}{B_2}:=\frac{A_3}{B_3} \pointloin$$

\begin{proposition}
\label{prop-def_plus}
Let $\F=\udisj\limits_{n\geq 0} \F_n$ be the scheme of pointed rational functions. Then the above morphisms define a graded monoid structure on $\F$:
  $$ \plus: \F \times \F  \fd \F \pointloin$$
\end{proposition}

\begin{proof}
The only thing to prove is that the associativity condition
is satisfied. It is a consequence of the associativity
 of matrix multiplication. 
\end{proof}

\begin{remarques} 
\label{rem_loi_plus_sur_pi0}
\begin{enumerate} 
\item The above monoid structure on $\F$ induces a graded monoid structure on its connected components 
$(\pN\F)(k):=\udisj\limits_{n\geq 0}(\pN\F_{n})(k)$, and thus on $\ppn$. The monoid law  on these sets is again
 denoted by $\plus$.
\item \label{homo_non_triviale}Taking the $\plus$-sum of ``trivial'' homotopies of rational functions
produces ``non trivial'' homotopies.  
\item \emph{Warning}: we use additive notations for the monoid law on 
$\F$, but we would like to stress that it
is \emph{non commutative}. However, we will see 
in corollary~\ref{corollaire_plus_pi0_commut} that
the monoid  $(\ppn, \plus)$ is abelian. 
\end{enumerate}
\end{remarques}

\begin{exemple} 
\label{ex_fc}
We now give some examples of  $\plus$-sums of rational functions.
\begin{enumerate}
\item One has 
$$ {X} \plus {X} = \frac{X^2-1}{X}$$
\item For any pointed rational function
$\frac{A}{B}$, one has
$$ {X} \plus \frac{A}{B} = \frac{AX -B}{A}  \et 
\frac{A}{B} \plus \frac{X}{1} = \frac{AX-V}{BX+U} \pointloin$$
\item \label{fc} More generally, for every monic polynomial $P\in k[X]$, and for every unit
$b_0 \in k^{\times}$, one has
$$ \frac{P}{b_0} \plus \frac{A}{B} = \frac{AP -\frac{B}{b_0}}{b_0A} =\frac{P}{b_0}- \frac{1}{b_0^2\frac{A}{B}} \qquad.$$ 
\end{enumerate}
\end{exemple}

\begin{remarque} 
\label{rem_fractions_continues}
\label{tautologique}
When dealing with rational functions with coefficients in the field
$k$, here is a convenient way to understand  the preceding addition law $\plus$ on the level of continued fraction expansions.
Every pointed rational fraction $f=\frac{A}{B} \in \F_n(k)$
 admits a unique \emph{twisted} continued fraction expansion 
of the following form:
$$ \frac{A}{B}= \frac{P_0}{b_0} - \frac{1}{b_0^2\left(\frac{P_1}{b_1}- \frac{1}{b_1^2(\dots)} \right)} \qquad, $$
where for each $i$, $P_i\in k[X]$ 
is a monic polynomial of positive degree and $b_i$ is a non-zero scalar in $k^\times$.
Such an expansion  always stops, as the sum of the degrees of the $P_i$ equals the degree of $A$.\\
 Example~\ref{ex_fc} then shows that 
the twisted continued fraction expansion of the $\plus$-sum of
two rational functions is the concatenation of the two twisted
continued fraction expansions.\\
For the sequel, it is useful  to note that every 
 rational function $f\in \F_n(k)$
is tautologically the $\plus$-sum of the polynomials appearing in 
its twisted continued fraction expansion:
$$f=\frac{P_0}{b_0} \plus \frac{P_1}{b_1} \plus \dots \plus \frac{P_r}{b_r} \qquad.$$
\end{remarque}

\subsection{The Bézout form}
\label{subsec_bezout}

In the $18^{\mathrm{th}}$ century, Bézout described a way to associate to every
rational function a non-degenerate symmetric matrix\footnote{
 It is sometimes referred as ``the bezoutian'' of the rational function,
 but we prefer the terminology ``Bézout form''.}. In modern terms, 
Bézout's construction describes, for every positive integer $n$,
 a scheme morphism  
$$ \bez_n : \F_n \fd \mc{S}_n \qquad,$$
where $\mathcal{S}_{n}$
is the scheme of non-degenerate \nn{} symmetric matrices.
The homotopy invariants 
$ (\pN \F_n)(k)
\xrightarrow{ \pN \bez_n}
  (\pN \mc{S}_n)(k)$
associated to these morphisms (see remark~\ref{rem_invariants})
are at the center of our study.\\ 

We start first by reviewing Bézout's construction.

\begin{definition}
\label{def_forme_bezout}
Let $R$ be a ring, $n$ be a positive integer and 
$f=\frac{A}{B}$ be an element of $\F_n(R)$.\\
 The polynomial  $A(X)B(Y)-A(Y)B(X) \in R[X,Y]$ is
 divisible by $X-Y$. 
Let  
$$\delta_{A,B}(X,Y):=\frac{A(X)B(Y)-A(Y)B(X)}{X-Y}=:\sum\limits_{1\leq p, q\leq n} c_{p, q}X^{p-1} Y^{q-1} \pointloin$$
Observe that the coefficients of $\delta_{A,B}(X,Y)$ are
symmetric in the sense that one has 
$$ c_{p,q}=c_{q,p} \qquad \forall \ 1\leq p, q\leq n.$$
The Bézout form of $f$ is the symmetric bilinear form over $R^n$
whose Gram matrix is the \nn-symmetric matrix $\left[c_{p,q}\right]_{1\leq 
p,q\leq n}$. We denote it $\bez_n(A,B)$ or $\bez_n(f)$ and  
Bézout's formula: 
\begin{equation}
\label{bez_formula}
\det \bez_n(f) = (-1)^{\frac{n(n-1)}{2}} \res(f)
\end{equation}
shows that this form is non-degenerate.\\
The above construction describes for every positive integer a natural transformation of functors
$ \F_n(-) \fd \mc{S}_n(-)$
and thus a  morphism of schemes
$$ \bez_n:\F_n \fd \mc{S}_n \pointloin$$
\end{definition}

\begin{remarque}
Following \cite{GZK}, chapter~III, example~4.8, we give 
now a  more conceptual definition of the Bézout form 
in terms of Serre duality.\\
Let $\frac{A}{B}$ be a pointed degree $n$ 
rational function.
Consider the complex of vector bundles over $\pun$: 
$$\mathcal{O}(-n-1) \xrightarrow{\begin{bmatrix} 
                                  -B\\A
                                 \end{bmatrix}
} \mathcal{O}(-1) \oplus \mathcal{O}(-1)
\xrightarrow{\begin{bmatrix} 
                                  A &B 
                                \end{bmatrix}} \mathcal{O}(n-1) \pointloin$$ 
The associated hypercohomology spectral sequence has only one non-trivial
differential, $\mathrm{d}_2^{0,1}: \H^1(\pun, \mc{O}(-n-1)) \fd 
\H^0(\pun, \mc{O}(n-1))$. The composite
$$ \bp{\H^0(\pun, \mc{O}(n-1))}^* \isom \H^1(\pun, \mc{O}(-n-1)) \fd 
\H^0(\pun, \mc{O}(n-1))$$ 
identifies with the Bézout form.
\end{remarque}

\subsection{The main theorem}
 We can now state our main result.
\label{subsec_enonce}

\begin{theoreme}
\label{th_fr}
 The Bézout invariants  distinguish all the naive homotopy classes
of rational functions: the map
$$   \Bp{ \udisj\limits_{n\geq 0} (\pN \F_n)(k) , \plus }
\xrightarrow{\udisj\limits_{n\geq0} \pN \bez_n}
 \Bp{\udisj\limits_{n\geq 0} (\pN \mc{S}_n)(k), \oplus }$$
is  an isomorphism of graded monoids.
\end{theoreme}

\begin{corollaire}
\label{corollaire_plus_pi0_commut}
Since the monoid $\bp{\udisj\limits_{n \geq 0}(\pN\mc{S}_n)(k), \oplus}$ is abelian,  
$\bp{\ppn, \plus}$ is also abelian.
\end{corollaire}

In order to make the theorem more explicit, we now describe the 
set of values $(\pN \mc{S}_n)(k)$ of the Bézout invariants.
To do so, we need the following definition.

\begin{definition} 
\begin{enumerate}
\item The Witt monoid of the field $k$ is the monoid, 
for the orthogonal sum $\oplus$, 
of isomorphism classes of non-degenerate symmetric $k$-bilinear forms.
We denote it $\MW(k)$.
\item Let $\tMW(k)$ be the monoid of 
\emph{stable}\footnote{If the field $k$ is not of characteristic $2$, $\tMW(k)
$ is equal to ${\MW}(k)$. So a reader not interested in
 the case of a base field of characteristic $2$ can forget about stabilization.}
 isomorphism
classes of non-degenerate symmetric $k$-bilinear forms. By definition, 
it is the quotient of $\MW(k)$ where two forms $b$ and $b^\prime$
 are to be identified if there exists a form $b^{\prime\prime}$ such that
$b \oplus b^{\prime \prime} \isom b^{\prime} \oplus b^{\prime \prime}$.
It comes with a natural grading induced by the rank,
and for every positive integer $n$, we denote $\tMW_n(k)$ the degree $n$-component of $\tMW(k)$. 
\item  The Grothendieck-Witt group $\GW(k)$ is 
 the Grothendieck group of the monoid $\MW(k)$ (or of $\tMW(k)$).
\end{enumerate}
 \end{definition}

The description of the sets $(\pN \mc{S}_n)(k)$ is given next.
\begin{proposition}
\label{prop_pi0Sn}
Let $n$ be a positive integer.
 \begin{enumerate}
\item \label{invariant_dans_GW} The canonical quotient map
$q_n: \mc{S}_n(k) \fd \tMW_n(k)$
factors through $(\pN \mc{S}_n)(k)$:
$$\xymatrix{
\mc{S}_n(k) \ar[r]^-{q_n} \ar[d]& \tMW_n(k) \\
(\pN \mc{S}_n)(k) \ar[ur]_-{\overline{q}_n}
}$$
\item \label{pi0Sn}
Let
$$\tMW_n(k) \croix{\kk} k^\times$$
be the canonical fibre product induced by the \emph{discriminant} map $\tMW_n(k)\fd \kk$. 
Then the  map 
$$ \Bp{  \udisj\limits_{n\geq0}(\pN\mc{S}_n)(k), \oplus }  \xfd{\udisj\limits_{n\geq0}\overline{q}_n \times \det} \Bp{ \udisj\limits_{n\geq0}\tMW_n(k) \croix{\kk} k^\times, \oplus}  \pointloin$$
is a monoid isomorphism. (Above, the right-hand term is endowed with the
canonical monoid structure induced by the orthogonal sum in $\tMW(k)$ and the product in  $k^\times$).
\end{enumerate}
\end{proposition}

Proposition~\ref{prop_pi0Sn} above is certainly well-known to specialists.
 An elementary proof, due to Ojanguren \cite{ojan},  
 is postponed until appendix~\ref{sec_hermite}, because it would be too digressive here.
The conjunction of
theorem~\ref{th_fr} and proposition~\ref{prop_pi0Sn} gives the 
following more explicit description of $\ppn$.

\begin{corollaire}
\label{coro_th_effectif}
There is a canonical isomorphism of graded monoids: 
$$ \left(\ppn, \plus\right) \isom   \Bp{\udisj\limits_{n\geq0}\tMW_n(k) \croix{\kk} k^\times, \oplus} \pointloin$$
\end{corollaire}

\begin{exemple} 
\begin{enumerate}
\item When $k$ is algebraically closed, we have an isomorphism of monoids
 $$
  \ppn \xrightarrow[\deg  \times \ \res]{\isom} \mathbf{N} \times k^\times \qquad.$$
  The  Bézout invariant is just the conjunction of the resultant and degree invariants.
\item When $k$ is the field of real numbers $\R$, we have an isomorphism of monoids:
 $$\ppn \xrightarrow[(\mathrm{sign}\, \rond \, \bez) \times \ \res]{\isom} (\mathbf{N} \times \mathbf{N})  \times \R^{\times} \qquad,$$
 $\mathrm{sign}$ denoting the signature of a real symmetric bilinear form.
In this case, the Bézout invariant is sharper than the resultant and the degree invariants.
  \end{enumerate}
\end{exemple}

\subsection{Proof of theorem~\ref{th_fr}}
\label{subsec_dem}
Three things are to be proved: injectivity, surjectivity  and 
compatibility with the monoid structures of the Bézout map. Each one has an elementary proof. 
In \S\ref{surj}, we start by proving the surjectivity and the compatibility condition,
both at the same time. 
We then remark in \S\ref{inj} that injectivity reduces
 to that of $\bez_2$. We conclude in \S\ref{n2} by 
an independent analysis of the scheme $\F_2$ and of the map $\bez_2$.

\subsubsection{Surjectivity and monoidal compatibility of the Bézout map}
\label{surj}

It is very convenient in the sequel to adopt the following conventions.
\begin{definition}
 \label{notation_fraction_et_forme}
Let $n$ be a positive integer. For every list of 
units $u_{1}, \dots, u_{n} \in k^\times$, let
\begin{itemize}
\item $\langle u_1, \dots, u_n \rangle$
denote the diagonal symmetric bilinear form 
$\langle u_1 \rangle \oplus \dots \oplus \langle u_n \rangle \in \mc{S}_n(k)$.
\item $[u_1, \dots, u_n]$ denote the pointed rational function 
$ \frac{X}{u_1} \plus \dots \plus \frac{X}{u_n} \in\F_{n}(k)$.
\end{itemize}
\end{definition}
 
The following lemma shows that, up to naive homotopy, 
any symmetric bilinear form and any rational function is of
 the preceding form.

\begin{lemme}
\label{lemme_somme_monomes}
Let $n$ be a positive integer. Then:
\begin{enumerate}
\item For any symmetric bilinear form  $S\in \mc{S}_n(k)$ 
there exists units $u_1, \dots ,
u_n \in k^\times$ such that $S$ is homotopic to the diagonal form $ \langle u_1,\dots,u_n \rangle$.
 \item For any pointed rational function $f\in \F_n(k)$ 
there exists units $u_1, \dots ,
u_n \in k^\times$ such that we have
$$f\hp [u_1,\dots,u_n] \qquad .$$
\end{enumerate}
\end{lemme}

\begin{proof}\begin{enumerate}
\item Suppose first that the  characteristic of the field $k$ is not $2$. Then  every symmetric matrix
$S \in \mc{S}_n(k)$ is conjugate by an element $P \in \SL_n(k)$ to a 
diagonal matrix. Decomposing $P$ into a product of elementary matrices 
and multiplying by $1-T$ the non diagonal elements of these elementary matrices yields a homotopy (that is to say an element of $\mc{S}_n(k[T])$)
to a diagonal matrix.\\
If $k$ is of characteristic $2$, then $S$ is conjugate by an element 
$P \in \SL_n(k)$ to a block diagonal matrix, with possible $\begin{bmatrix} 
0 & 1 \\
1 & 0 
\end{bmatrix}$ terms. In addition to the preceding argument, one can use the homotopy $\begin{bmatrix} 
T & 1 \\
1 & 0 
\end{bmatrix}$ to link $S$ to  a diagonal matrix.
\item We prove this point by induction on the degree $n$ of $f$.\\
As noticed in remark~\ref{rem_fractions_continues}, a rational function 
$f \in \F_n(k)$ 
is tautologically the $\plus$-sum of some polynomials, say
$P_1 \plus \dots \plus P_k$.  
Thus one can assume that $f$ is a polynomial.
Example~\ref{exemple_homotopies} (\ref{exemple_homotopies_polynome})
then shows that a polynomial is always homotopic to its leading 
term. So it's enough to treat the case of a monomial
$\frac{X^n}{u}$, with $u\in k^\times$. Now,
example~\ref{exemple_homotopies}~(\ref{exemple_homotopies_inverse_polynome})
 shows that the element $\frac{X^n}{TX^{n-1}+u} \in \F_n(k[T])$ 
defines a homotopy between  $\frac{X^n}{u}$
and $\frac{X^n}{X^{n-1}+u}$.
But this last rational function decomposes as 
$X \plus g$ for some $g\in \F_{n-1}(k)$.
One concludes by using the inductive hypothesis on $g$. \end{enumerate}
\end{proof}

The monoids $((\pN \mc{S})(k), \oplus)$ and
$((\pN \F)(k), \plus)$ are generated by their degree $1$ components.
Since the map $\pN \bez_1: (\pN \F_1)(k) \fd (\pN \mc{S}_1)(k)$ is a bijection,
the following lemma shows that 
the Bézout form of a rational function of the form
  $[u_1,\dots,u_{n}]\in  \F_n(k)$ is
 homotopic to the diagonal form $\langle u_{1},\dots, u_{n}\rangle$.

\begin{lemme}
\label{bez} \label{bez_de_monoides}
Let $\frac{A}{B}\in \F_n(k)$ and $u\in k^\times$.
Then the Bézout form of $\frac{X}{u} \plus \frac{A}{B}$
 is conjugate\footnote{And is thus also homotopic!} by an element in $\SL_{n+1}(k)$ 
to the block diagonal form $\langle u \rangle \oplus \bez_n(A,B) $.
\end{lemme}

\begin{proof}
By definition, one has $\frac{X}{u} \plus \frac{A}{B}= \frac{XA-\frac{B}{u}}{uA}$.
Using the notations introduced in 
definition~\ref{def_forme_bezout}, we have
$$\delta_{XA-\frac{B}{u},uA}(X,Y)=u A(X)A(Y) +\delta_{A,B}(X,Y) \pointloin$$
In the basis  $(1,X,\dots, X^{n-1}, A(X))$, 
the matrix of the Bézout form is diagonal as announced.
\end{proof}
This proves that the Bézout map induces a surjective morphism of monoids. 
\subsubsection{Injectivity}
\label{inj}
Let $n$ be a positive integer. To prove the injectivity of
the map $\pN \bez_n$, we prove the injectivity of the composite
\emph{$$ (\pN \F_n)(k) \xrightarrow{\pN \bez_n} 
(\pN \mc{S}_n)(k) \xrightarrow{ \overline{q}_n \times \det}
 \tMW_n(k) \croix{\kk} k^\times \pointloin$$}
Because up to homotopy any rational function is a
 $\plus$-sum of degree 1 monomials (\cf lemma~\ref{lemme_somme_monomes}), the injectivity of the previous map can be reformulated as follows.

\begin{proposition}
\label{prop_inj}
Let $u_1, \dots, u_n, v_1, \dots, v_n$ be a list of units in  $k^\times$.
If the classes in $\tMW_n(k) \croix{\kk} k^{\times}$ of the diagonal forms  $\langle u_1,\dots,u_n\rangle$ and $\langle v_1,\dots,v_n\rangle$ are equal, then 
 $[u_1,\dots,u_n]\hp[v_1,\dots,v_ n]$ holds in $\F_n(k)$.
\end{proposition}

\begin{proof}
Let's first introduce some \emph{ad hoc} terminology.\\
Two diagonal forms
$\langle u_1, \dots, u_n \rangle$  and $\langle v_1, \dots, 
v_n \rangle$ are said to be equivalent through an 
\emph{elementary $\SL_2(k)$-transformation} if there
 exists an subscript  $1\leq i \leq n-1$ such that:
\begin{itemize}
\item the $2$-forms $\langle u_i, u_{i+1}\rangle$ and $\langle v_i, v_{i+1} \rangle$ are $\SL_2(k)$-equivalent;
\item for all  $j \neq i, i+1$, we have
$u_j=v_j$. 
\end{itemize}
\smallskip

The next lemma is a slight reformulation of \cite{milnor}, 
chapter III, lemma~5.6 which is gives 
a presentation of the Witt group $\mathrm{W}(k)$ by 
generators and relations.

\begin{lemme}
\label{elem}
Let $\langle u_{1},\dots, u_{n} \rangle$ and $\langle v_{1},\dots, v_{n}\rangle$ 
be two  diagonal forms in $\mc{S}_n(k)$.
Their images in $\tMW_n(k)\croix{\kk}k^{\times}$ are equal
if and only if one can ``pass from one to the other'' by 
a finite sequence of elementary $\SL_2(k)$-transformations.
\end{lemme}

So, in order to prove proposition~\ref{prop_inj}, 
it's enough to deal with the case $n=2$. 
This case is analysed independently in the next paragraph.

\subsubsection{Rational functions of degree $2$}
\label{n2}

Let $\ga$ be the ``additive group'',
 that is to say the affine line $\aun$ seen as a group scheme.
For every positive integer $n$, $\ga$ acts freely on 
$\F_n$ by translations, \ie{} by the formula
  $$ h \cdot \frac{A}{B}:= \frac{A+hB}{B} \qquad,$$
on the level of points.
The following lemma shows that the $\ga$-torsor 
$\F_n$ splits.

\begin{lemme}
 Let $R$ be a ring and $\frac{A}{B}$ be an element of $\F_n(R)$.
There exists a unique pair of polynomials $(U_1,V_1)$ of $R[X]$
with $\deg(U_1) = n-1$, $\deg(V_1) \leq n-1$ and such that 
$AU_1 + BV_1= X^{2n-1}$. Let 
$\phi_n\left(\frac{A}{B}\right)$ be the opposite of the coefficient of $X^{n-1}$ in $V_1$. Then the associated scheme morphism 
$$\phi_n: \F_n \fd \aun$$
 is $\ga$-equivariant.
 In particular, $\F_n$ splits as the product $\phi_n^{-1}(0) \times \aun$.
\end{lemme}

\begin{proof}
Let $A,B,U_1$ and $V_1$ be like above. If one changes $(A,B)$ to
$(A+hB,B)$, then $(U_1,V_1)$ is changed to
$(U_1, V_1-hU_1)$. The claim follows since $U_1$ is necessarily monic.
\end{proof}

Moreover, observe that the morphism $\bez_n: \F_n \fd \mc{S}_n$ is by construction $\ga$-equivariant
when $\mc{S}_n$ is endowed with the trivial action.
In dimension $2$, the morphism $\bez_2$ induces an isomorphism between $\phi_2^{-1}(0)$ and $\mc{S}_2$.

\begin{proposition}
\label{prop_F2}
 The morphism
$$ \F_2 \xrightarrow{\bez_2 \times \phi_2} \mc{S}_2 \times \aun$$
is a $\ga$-equivariant isomorphism of schemes.
\end{proposition}

\begin{corollaire}
\label{coro_Bez2_inj}
 The map
$$(\pN \F_2)(k) \xrightarrow{\pN \bez_2}(\pN \mc{S}_2)(k)
\xrightarrow{\overline{q}_2 \times \det } \tMW_2(k)\croix{\kk}k^{\times}  $$
is injective.
\end{corollaire}

Corollary~\ref{coro_Bez2_inj} concludes the proof of proposition~\ref{prop_inj}.\\

\begin{proof}[Proof of proposition~\ref{prop_F2}]
One can write down the inverse morphism $\psi:\mc{S}_2 \fd \phi_2^{-1}(0)$ by solving a system of two equations with  
two unknowns. The formula is
$$ {\tiny \psi \left(\begin{bmatrix}
\ \alpha & \beta\ \\
\ \beta & \gamma \
\end{bmatrix}\right)}
=
\frac{X^{2}+\frac{\alpha\beta}{\beta^{2}-\alpha\gamma} X+\frac{\alpha^2}{\beta^{2}-\alpha\gamma}}{\gamma X+\beta} \pointloin$$
\end{proof}

This proof of theorem~\ref{th_fr} is now complete.\end{proof}

\begin{remarque}
 For $n\geq 3$, proposition~\ref{prop_F2}  generalizes to an isomorphism
$\F_n \isom \mc{H}_n \times \aun$, where $\mc{H}_n$
is the scheme of non-degenerate \nn{} \emph{Hankel matrices}, 
see appendix~\ref{sec_hankel}. The miracle that allows our computation 
is that non-degenerate $2 \times 2$ 
Hankel matrices identify with regular non-degenerate symmetric matrices 
in which naive homotopies are well understood. 
\end{remarque}

\subsection{Comparing  naive and  motivic homotopy classes}
\label{subsec_comparaison}

Recall that  $\ho(\esp)$ is the homotopy category
 of spaces of $\aun$-homotopy theory defined by Morel and
 Voevodsky in \cite{MV}.
For any space $X \in \esp$, let 
 $\susp^sX$ be its ``suspension  with respect to the simplicial circle'', that is to say
$$ \susp^s X := \hocolim{} \left( \vcenter{\xymatrix@R=1ex{
{}\pt & & {}\pt  \\
& X \ar[ul] \ar[ur] &
}}\right) \pointloin$$

\begin{lemme}
In the homotopy category $\ho(\esp)$, there is an equivalence
$$ \pun \approx \susp^s \gm \qquad.$$
 Above, $\gm$ is the multiplicative group, that is to say the
 scheme $\aun-\{0\}$.
\end{lemme}

\begin{proof}
This is a consequence of $\pun$ being covered by two
contractible open subschemes (two affine lines $\aun$)
intersecting along a $\gm$. 
\end{proof}

It follows formally from the preceding lemma that 
the set $\pp^\mot$ is equipped with a 
\emph{group structure}, whose law is denoted by $\oplus^\mot$.
The canonical map $\ppn \fd \pp^\mot$ is thus not a bijection: 
its reaches only elements of non-negative degrees.
However,  Morel's computation shows that the error in the naive computation is
as small as possible.

\begin{theoreme}[Morel, \cite{morel}, theorem~4.36]
\label{th_morel}
There is a group isomorphism
$$ \bp{\pp^\mot,\oplus^\mot} \isom \bp{\GW(k) \croix{\kk} k{^\times}, \oplus} \qquad.$$
In particular, the group $\bp{\pp^\mot,\oplus^\mot}$ is (abstractly) isomorphic to 
the group completion of the monoid $(\ppn,\plus)$.
\end{theoreme}

As one can expect, the group completion is induced by the
canonical map between the monoid of naive homotopy classes to the group of $\aun$ homotopy classes. 

\begin{theoreme}
\label{conjecture_comparaison}
The canonical map 
$$ \bp{\ppn, \plus} \fd \bp{\pp^\mot, \oplus^\mot }$$
 is a group completion.
\end{theoreme}

\begin{proof}
The main point is to prove the compatibility
between the monoid law $\plus$ on $\ppn$ and the group
law on $\pp^\mot$. Namely, we have the following
\begin{proposition}
\label{prop:plus_oplus_texte}
 The canonical map
$$ \bp{\ppn , \plus} \fd \bp{\pp^\mot, \oplus^\mot}$$
is a monoid morphism.
\end{proposition}
A detailed proof of proposition~\ref{prop:plus_oplus_texte} is 
postponed to appendix~\ref{appendice_comparaison}. 

\end{proof}

\begin{remarque}
The question of the comparison between
 the naive addition law $\plus$ and the $\aun$-addition law
$\oplus^\mot$ makes sense on the level of the spaces of
 morphisms and not only on the level of homotopy classes. 
In particular, it is natural to ask
if for any pair of non-negative integers $(m,n)$ 
the following diagram homotopy commutes\footnote{
{In the diagram, $\Omega^{\pun} \pun$ is the space (in $\esp$) of 
pointed morphisms from $\pun$ to itself. This space splits
as the disjoint union of its components of constant degree: 
$\Omega^{\pun} \pun= \udisj\limits_{k\in \Z}
\Omega^{\pun}_k \pun$. For any non-negative $k$, $\F_k$ 
is seen as a subspace of $\Omega^{\pun}_k \pun$.}
}:
$$\vcenter{\xymatrix{
\F_m \times \F_n \ar[r]^-{\plus} \ar@{^(->}[d] & \F_{m+n} \ar@{^(->}[d] \\
\Omega^{\pun}_m \pun \times \Omega^{\pun}_n \pun \ar[r]^-{\oplus^\mot} & \Omega^{\pun}_{m+n} \pun
}} \pointloin$$ 
\end{remarque}

\section{Related computations}
\label{sec_extensions}

The previous computation has natural extensions which we
consider now.
\begin{itemize}
 \item In \S\ref{libre}, we compute  the set of free
 homotopy classes of scheme endomorphisms of $\pun$.
\item In \S\ref{composition}, we 
make explicit  the monoid structure 
induced on $\ppn$ 
by the composition of endomorphisms of $\pun$.
\item Finally, in \S\ref{p1pd}, we compute the set of naive homotopy classes of morphisms from $\pun$ to $\pd$.
\end{itemize}

\subsection{Free homotopy classes of rational functions}
\label{libre} \label{sec_libre} \label{subsec_libre}
We compute here the set $\ppl$ of naive homotopy classes
 of \emph{unpointed} endomorphisms of $\pun$. 
The result is mostly a consequence of our previous 
computation of $\ppn$.\\

The description of free endomorphisms of $\pun$ in terms of 
rational functions is the following.

\begin{definition}
 For every non-negative integer $n$, the scheme $\libre_n$ of unpointed degree $n$ rational functions  is the
open subscheme of $\mathbf{P}^{2n+1}:= \proj (k[a_0,\dots, a_n, b_0, b_n])$ complementary to the hypersurface of equation
$$\res_{n,n}(a_nX^n + \dots + a_0, b_nX^n+ \dots + b_0)=0 \pointloin$$
\end{definition}

\begin{proposition}
\label{prop-def_libre}
\begin{enumerate}
\item The datum of an endomorphism of $\pun$ is equivalent to the
datum of a non-negative integer $n$ (its \emph{degree}) and
of an element in $\libre_n(k)$.
\item Moreover, the datum of a naive homotopy $\pun \times \aun \fd \pun$ is equivalent to the datum of a non-negative integer $n$ and of an element in $\libre_n(k[T])$. 
\end{enumerate}
\end{proposition}

We denote $\hl$ the equivalence relation generated
by \emph{unpointed} naive homotopies and by $\ppl$ the 
set of unpointed naive homotopy classes of endomorphisms of $\pun$.
It follows from proposition~\ref{prop-def_libre}
that the degree is still a naive homotopy invariant.
Thus the set $\ppl$ splits as the disjoint union
$$\ppl=\udisj\limits_{n\geq0} \ppl_n$$
of its components of fixed degree.\\

Let $n$ be a positive integer. Note that a $k$-point in 
$\libre_n$ gives a pair of coprime polynomials $(A,B)$
only up to the multiplication by a unit of $k$. Thus the Bézout 
form of an unpointed endomorphism of $\pun$ is defined only
up to multiplication by the square of a unit of $k$.
The result is that the Bézout form leads again to
an invariant that distinguishes all the homotopy classes.

 \begin{theoreme}
\label{theoreme_libre}
The canonical map of graded sets:
$$
   \begin{array}{rccl}
        \ppl=& \udisj\limits_{n \geq 0} \ppl_n & \fd &  \udisj\limits_{n \geq 0} \tMW_n(k) \croix{\kk} \kkn \\
       {}& \left[\frac{A}{B}\right] &\mapsto & \left[ \bez(A,B),\ \det \bez(A,B) \right]
\end{array}
$$
is a bijection.
 \end{theoreme}
 
 \begin{proof}[Proof of theorem~\ref{theoreme_libre}]
 This is a consequence of corollary~\ref{coro_th_effectif} and
of the following lemma.
 \end{proof}

 \begin{lemme} 
\label{lemme_libre}
 \begin{enumerate}
 \item Any free rational function is naively homotopic to a 
pointed rational function.
 \item Let $f$ and $g$ be two \emph{pointed} rational functions. Then one has the relation  $f \hl g$ 
if and only if there exists a non-zero element 
 $\lambda \in k^\times$ with  $f \hp \lambda^2 g$.
 \end{enumerate}
 \end{lemme}
 
 \begin{proof} 
 \begin{enumerate}
 \item Let $f=\frac{A}{B}$ represent a rational function
and let  $\alpha_1$ be a matrix in $\SL_2(k)$ such that 
$\alpha_1 \cdot \infty = f(\infty)$, for the usual action
of $\SL_2(k)$ on $\pun(k)$.
 Let $\alpha(T)$ be an algebraic path in 
 $\SL_2(k[T])$ linking the identity to $\alpha_1$. (Again, this
can be obtained using a decomposition of $\alpha_1$ as a product 
of elementary matrices). The column vector
 $\alpha(T)^{-1}\cdot \begin{pmatrix} A \\ B \end{pmatrix}$ is a $k[T]$-point of $\libre_n$ and thus 
yields a homotopy between $f$ and the \emph{pointed} rational function  $\alpha_1^{-1} \cdot \frac{A}{B}$.
 \item Suppose first that we have two rational functions $f,g$ such that  $f \hl g$ and let's show that there exists a unit $\lambda \in k^\times$ such that
 $f \hp \lambda^2 g$.\\
 Let
 $$f=f_0 \underset{F_0(T)}{\hl} f_1 \underset{F_1(T)}{\hl} \dots \underset{F_{N-1}(T)}{\hl} f_N=g$$
 be a sequence of elementary homotopies between $f$ 
and $g$. 
Let  $\alpha(T)$ be an matrix in $\SL_2(k[T])$
such that we have $F_0(T,\infty)= \alpha(T) \cdot \infty$
and $\alpha(0)=\id$. The path 
 $\alpha(T)^{-1} \cdot F_0(T)$ 
yields a pointed homotopy between the pointed rational 
functions  $f_0$ and $\alpha(1)^{-1} \cdot f_1$.
Moreover, for $N>1$, we then have a  
sequence of $N-1$ free homotopies 
 $$\alpha(1)^{-1}\cdot f_1 \underset{\alpha(1-T)^{-1} F_1(T)}{\sim} f_2 \underset{F_2(T)}{\sim} \dots \underset{F_{N-1}(T)}{\sim}f_N=g \pointloin$$
Thus the result will follow by induction from the 
case $N=1$.\\
 When $N=1$, $\alpha(1)^{-1}$ is of the form
 $\begin{bmatrix} \lambda & \mu \\ 0 & \lambda^{-1} \end{bmatrix}$.
 So, 
 $\begin{bmatrix} 1 & T \lambda \mu \\
 0 & 1\end{bmatrix} \cdot \alpha(T)^{-1} \cdot F_0(T)$ 
gives a pointed homotopy between $f_0$ and $\lambda^2 f_1$.\\
We now show the converse. For every $\lambda \in k^\times$,
a path in $\SL_2(k[T])$ between the identity matrix and 
 $\begin{bmatrix} \lambda & 0 \\ 0 & \lambda^{-1} \end{bmatrix}$ 
yields a free homotopy between any rational function $g$ and  $\lambda^2 g$. The result follows.
 \end{enumerate}
 \end{proof}

\subsection{Composition of rational functions} 

\label{composition} \label{subsec_composition}

The set $\ppn$ 
admits a second monoid structure induced by composition of
morphisms. We make explicit this new structure in terms of our
previous description of $\ppn$.

\begin{definition}
Define a new composition law, say $\rond$, on 
 \emph{$\udisj\limits_{n\geq 0}\tMW_n(k) \croix{\kk} k^\times$}
 by
$$ (b_1, \lambda_1) \rond (b_2,\lambda_2) :=
\big{(}b_1 \tens b_2, \lambda_1^{\dim b_2} \lambda_2^{(\dim b_1)^2}\big{)} \qquad.$$
This law is associative and endows
 \emph{$\udisj\limits_{n\geq 0} \tMW_n(k) \croix{\kk} k^\times$} with a monoid structure.
\end{definition}

\begin{theoreme}
\label{th_composition}
The following map induces an isomorphism of graded 
``bi-monoids''
\appsn{\left(\ppn, \plus, \rond \right)}
{ \Big{(}\udisj\limits_{n\geq 0}\tMW_n(k) \croix{\kk} k^\times, \oplus, \rond  \Big{)}}
{f}{[\bez(f), \res(f)]}
\end{theoreme} 

\begin{remarque} \attention the triple
 $(\ppn, \plus, \rond)$
is \emph{not} a semi-ring. In general, one has
 distributivity of $\rond$ over $\plus$ only on the left-hand side. That is to say that for any triple
 $(f, g_1, g_2)$ of pointed rational functions, we have
 $$ f \rond(g_1 \plus g_2) \hp (f\rond g_1) \plus (f\rond g_2) \qquad,$$
but in general 
$$ (g_1 \plus g_2) \rond f \hp (g_1 \rond f) \plus (g_2 \rond f) $$
\emph{does not hold}.
\end{remarque}

\begin{proof}[Proof of theorem~\ref{th_composition}]
Since any pointed rational function is up to homotopy a
$\plus$-sum of degree $1$ rational functions, theorem~\ref{th_composition} follows from the following lemma.

\begin{lemme} Let  $a \in k^\times$ be a unit, $m$ and $n$ be positive integers, and 
$f=\frac{A}{B} \in \F_m(k)$ and $g= \frac{C}{D}\in F_n(k)$
be two pointed rational functions. Then
\begin{enumerate}
\item We have $\frac{X}{a} \rond f = \frac{1}{a}f$.
  \item In the stable Witt monoid $\tMW_{(m+1)n}(k)$, we have
  $$ \bez_{(m+1)n}\left((X \plus f) \rond g \right) = \bez_n(g) \plus \bez_{mn}(f\rond g) \qquad.$$
  \item We have
  $$ \det \bez_{(n+1)m}\left((X \plus f) \rond g \right) =
  \det \bez_n(g)^{2m+1} \cdot \det \bez_{nm}(f \rond g) \qquad.
  $$
\end{enumerate}
\end{lemme}
\begin{proof}
\begin{enumerate}
\item This is true by definition.
\item Let $\dfrac{\tilde{A}}{\tilde{B}}$ 
be the pointed rational function representing $f \rond g$.
By definition, one has:
$$ \tilde{A}(X)= \sum_{i=0}^m a_i C(X)^iD(X)^{m-i}  \et \tilde{B}(X)=\sum_{i=1}^m b_i C(X)^i D(X)^{m-i}  
 \qquad.$$
Since we have $X \plus \frac{A}{B}= X - \frac{B}{A}$, 
 the endomorphism
$(X \plus f) \rond g$ is represented by the 
rational function 
$$\frac{C}{D}-\frac{\tilde{B}}{\tilde{A}} = \frac{C \tilde{A} - D \tilde{B}}{D \tilde{A}} \qquad.$$
Moreover, using the notations of definition~\ref{def_forme_bezout}, we have the identity
$$ \delta_{{C \tilde{A} - D \tilde{B}},{D \tilde{A}}}(X,Y)= \tilde{A}(X)\tilde{A}(Y) \delta_{C,D}(X,Y) + D(X)D(Y) \delta_{\tilde{A},\tilde{B}}(X,Y) \qquad.$$
Observing the congruence
$\tilde{A} \equiv C^m \mod D$, we deduce 
$$ \res(\tilde{A},D) = \res(C^m,D)= \res(C,D)^m \in k^\times \qquad.$$
Thus the family $(D(X), X \cdot D(X), \dots, X^{mn-1} \cdot D(X), \tilde{A}(X),X\cdot\tilde{A}(X), \dots, X^{m-1} \tilde{A}(X) )$ gives a basis of the $k$-vector space of polynomial of degree $<(n+1)m-1$. And in this basis, the matrix of the form $\bez_{(n+1)m}((X \plus f) \rond g)$ is 
$$ \begin{bmatrix} \bez_{mn}(\tilde{A}, \tilde{B}) & 0 \\ 0 &\bez_{m}(C,D)\end{bmatrix} \qquad.$$ 
\item This point follows from the proof of the previous one. Indeed, we
 have just proved the following  matrix identity
$$ \bez_{(n+1)m}(C \tilde{A} - D\tilde{B},D \tilde{A}) = \transp \syl_{}(\tilde{A}, D) \begin{bmatrix} \bez_{mn}(\tilde{A}, \tilde{B}) & 0 \\ 0 &\bez_{m}(C,D)\end{bmatrix} \syl(\tilde{A},D) \qquad,$$
where $\syl$ is the Sylvester matrix (see \cite{Bourbaki}, \S6,~$\mathrm{n}^{\rond}$6,~IV). The result now follows
from the relation $\det \syl_{}(\tilde{A}, D) = \res(\tilde{A},D)= \res(C,D)^m$. 
\end{enumerate}\end{proof}
This completes the proof of theorem~\ref{th_composition}.
\end{proof}

\subsection{Naive homotopy classes of morphisms to higher dimensional projective spaces}
\label{p1pd}
\label{subsec_p1pd}

Let $d\geq 2$ be an integer. 
We compute now
naive homotopy classes
 of morphisms from $\pun$ to $\pd$ (over $\spec k$).
The group of $\aun$-homotopy classes of morphisms from $\pun$ to $\pd$ is also determined in \cite{morel}
and we can again compare our result to it. 
 Not surprisingly, the computation is
much easier than the previous one for $d=1$.\\
 
 For brevity, 
we treat only the case of pointed morphisms. 
The base point in $\pd$ is taken at 
$\infty:= [1:0:\dots : 0]$. Pointed scheme morphisms 
from $\pun$ to $\pd$ and naive homotopies between 
them admit the following concrete description, analogous to 
that of proposition~\ref{prop_description} and proposition~\ref{prop_decodage_homotopie}.

\begin{definition}
 For every non-negative integer $n$, let $\F_n^d(-)$ be the functor from the  category of $k$-algebras to the category of sets  which maps any $k$-algebra $R$ to the set of pairs $(A,B)$, where $A$ 
is a monic polynomial of $R[X]$ of degree $n$, and  $B:=(B_1,\dots, B_d)$ is a $d$-tuple of polynomials, each  of degree strictly less than $n$ and such that the ideal  generated by the family $\{A,B_1,\dots,B_d\}$ is $R[X]$.
 These functors $\F_n^d(-)$ are representable by smooth schemes, which we denote $\F_n^d$ again.
\end{definition}

\begin{proposition}  
\begin{enumerate}
\item The datum of a pointed morphism $\pun \fd \pd$ is equivalent to the datum of a non-negative integer $n$, called the \emph{the degree} of the morphism, and of an element in $\F_n^d(k)$.
\item The datum of a pointed homotopy $\pun \times \aun \fd \pd$ is equivalent to the datum of a non-negative integer $n$ and of an element in $\F_n^d(k[T])$.
\end{enumerate}
\end{proposition}

\begin{remarque}
\label{rem_AH}
Let $R$ be a ring and $A$ be a monic polynomial in $R[X]$.
The datum of a $B$ such that $(A,B)$ is an $R$-point of $\F_n^d$ is equivalent to the datum of an element  
$$(\overline{B_1},\dots, \overline{B_d}) \in(\A^d-\{0\})\Bp{\sur{R[X]}{(A)}} \pointloin$$
 This point of view leads to the definition of \emph{Atiyah and Hitchin schemes}, see \cite{th_AH, these}.
\end{remarque}

 As in  definition~\ref{def_hp}, we denote  $\hp$ 
the equivalence relation
generated by pointed naive homotopies of morphisms from $\pun$ to $\pd$ and by 
$\ppdn$ the corresponding set of naive homotopy classes.
Since the degree of  morphisms is invariant through naive homotopies, the set  $\ppdn$ splits as  the union  
of its degree components:
$$ \ppdn = \udisj_{n \geq 0} \ppdn_n \qquad.$$

We can now state the result.
\begin{theoreme}
\label{th_d2}
For every $d \geq 2$, the degree map $\deg$
$$ \ppdn \xrightarrow[\isom]{\deg} \mathbf{N}$$
is a bijection.
\end{theoreme}

\begin{proof} Fix a non-negative integer $n$. We are going 
to prove that the set 
 $(\pN\F_n^d)(k)$ contains only one element.
 More precisely, 
we are going to link any  
element $(A, B_1, \dots, B_d) \in \F_n^d(k)$ to  $(X^n,1, \dots, 1)$ by a 
sequence of pointed  naive homotopies.\\
Notice first, that it's enough to link  
 $(A, B_1, \dots, B_d)$ to 
$(A,1, \dots, 1)$ since the pointed homotopy 
$$\left((1-T)A +TX^n , 1, \dots, 1\right)$$
 then links  $(A,1, \dots, 1)$ to $(X^n,1, \dots, 1)$.\\
 To do so, decompose the polynomial $A$ as a product of irreducible 
factors:
$$ A= \prod\limits_{i=1}^r P_i^{r_i} \qquad.$$
As noticed in remark~\ref{rem_AH}, the set of $k$-points of $\F_n^d$
with first coordinate equal to $A$ is in bijection with
$(\mathbf{A}^d -\{ 0 \}) \Bp{\sur{k[X]}{(A)}}$. 
The Chinese remainder theorem  identifies this set and
$$ \prod_{i=1}^r (\mathbf{A}^d-\{0\})\Bp{\sur{k[X]}{(P_i^{r_i})} } \pointloin$$
This means that it's enough to treat the case when 
$A$ is a power of an irreducible polynomial, say $A=P^r$.
In this case, the ring $R:=\sur{k[X]}{(P^r)}$ is local. 
A $d$-tuple $(\overline{B_1}, \dots, \overline{B_d}) \in R^n$ represents an element of
$(\A^d- \{0\})\Bp{\sur{k[X]}{(P^r)}}$
if and only if one of the $\overline{B_i}$ is a unit of $R$. 
Up to reordering, we can assume that $\overline{B_1}$ is such a unit.
The  $k[T]$-point of $\F_n^d$
$$(A, B_1, (1-T) B_2 +T, \dots, (1-T)B_d +T)$$ 
then yields a pointed homotopy linking
$(A, B_1,B_2, \dots, B_d)$ to $(A, B_1,1, \dots, 1)$.
One concludes using the pointed homotopy $(A, (1-T)B_1+T,1, \dots, 1)$ from $(A,B_1, 1,\dots, 1)$ to $(A,1,1, \dots, 1)$. 
\end{proof}

\begin{remarque}   It follows from theorem~\ref{th_d2}
 that the set $\ppdn$ is \emph{a fortiori}
endowed with a monoid structure, pulled back from that of $(\mathbf{N},+)$.
 We denote again its law $\plus$.
The statement of theorem~\ref{th_d2} would be much nicer
if there were an \emph{a priori} monoid structure on $\ppdn$ such
 that the degree map induces a monoid isomorphism. 
The author suspects that there might even be a monoid 
structure on the scheme 
$\udisj\limits_{n \geq 0}\F_n^d$.
\end{remarque}

For the same reason as before, $\ppd^\mot$ is a group. Morel's result is the following.

\begin{theoreme}[Morel, \cite{morel}, Theorem~4.13]
For every integer $d\geq 2$, the degree map
$$\bp{\ppd^\mot,\oplus^\mot} \xrightarrow{\deg} (\Z,+) $$
induces a group isomorphism.
\end{theoreme}

Naive and $\aun$-homotopy classes
compare again very well.

\begin{corollaire}
For every integer $d\geq 2$, the canonical map
$$ \bp{\ppdn, \plus}\fd \bp{\ppd^\mot, \oplus^\mot}$$
is a group completion. 
\end{corollaire}

\backmatter
\appendix

\section{Hermite inequality}
\label{sec_hermite}

The goal of this appendix is to present a proof,
 due to Ojanguren \cite{ojan}, of proposition~\ref{prop_pi0Sn}
which gave a description of the naive homotopy classes of non-degenerate
 symmetric bilinear forms. 
Since the reference  \cite{ojan} is not 
so common in libraries, we want to include this material here.
 The proof is  based on an elegant use of 
 Hermite inequality for symmetric bilinear forms over 
the principal ring $k[T]$.\\
 
More precisely, proposition~\ref{prop_pi0Sn} is a consequence of 
the following proposition, which gives a concrete
description of naive homotopies of non-degenerate symmetric bilinear forms.

\begin{proposition}
\label{pi0sl_revu}
Let $n$ be a positive integer and
let $S(T)$ be a \nn{} non-degenerate symmetric matrix
with coefficients in $k[T]$ (that is to say an element of $\mc{S}_n(k[T])$). Then, there exists a matrix  
$P(T) \in \SL_n(k[T])$ such that $\transp{P(T)} S(T) P(T)$
is block diagonal, where the block entries are whether units of $k$, whether 
$2 \times 2$-blocks of the form 
$\begin{bmatrix} 
0 & 1 \\
1 & \alpha(T)
\end{bmatrix}$
 with $\alpha(T) \in k[T]$.\\
In particular, since we have, for any $\alpha \in k$, the equality
$\begin{bmatrix} 
0 & 1 \\
1 & \alpha
\end{bmatrix} \equiv \begin{bmatrix} 
0 & 1 \\
1 & 0
\end{bmatrix}$ in $\tMW(k)$, the images of $S(0)$ and $S(1)$ in $\tMW(k)$
are equal.
\end{proposition}

\subsection{Hermite inequality}

We review first the classical version of Hermite inequality for integral
symmetric bilinear forms.

\begin{theoreme}[Hermite inequality]
\label{hermite_classique}
Let $n$ be a positive integer, 
 $L$ be a free $\Z$-module of rank $n$ and 
$b: L \times L \fd \R$
be a symmetric bilinear form (possibly degenerate).\\
Define the minimum of $b$ to be
$$\mu(b) := \min\limits_{x \in L - \{0\}} \vert b(x,x) \vert \qquad.$$
Then, one has the following inequality:
$$ \mu(b) \leq \left( \frac{4}{3}\right)^{\frac{n-1}{2}} \vert \D(b) \vert ^{\frac{1}{n}} \qquad, $$
where $D(b)$ stands for the \emph{discriminant}\footnote{
That is to say the determinant of the Gram matrix of $b$ in any base of $L$.
} of the form $b$.
\end{theoreme}

Hermite inequality is intimately related to the fact that 
any real number can be approximated by an integer at distance at most
$\frac{1}{2}$. It admits a generalization in the following algebraic context.

\begin{definition}
 Let $K$ be a field. An absolute value on $K$ is a function
\appsn{K}{\R^+}{x}{\vert x \vert}
satisfying the following conditions:
\begin{itemize}
\item $\vert x \vert =0 \iff x=0$;
\item $\forall x,y \in K, \quad \vert x y\vert = \vert x\vert  \vert y \vert$;
\item $\vert x+y \vert \leq \vert x \vert + \vert y \vert$. 
\end{itemize}
The absolute value is said to be 
\emph{non-archimedean} if we have furthermore
$$ \forall x,y \in K,\qquad \vert x+y \vert \leq 
\max(\vert x\vert, \vert y \vert) \pointloin$$ 
\end{definition}

\begin{prop-def}
\label{anneaudeHermite}
Let $K$ be a field with an absolute value. A subring
$R\subset K$  is of Hermite if it satisfies:
$$ \xymatrix@R=0ex{\forall a \neq 0, \quad \vert a \vert \geq 1\\
 \exists 0< \rho <1, \quad \forall x \in K, \quad
\exists r \in R, \quad  \vert x-r \vert \leq \rho}$$
The following  then holds:
\begin{itemize}
\item A Hermite subring $R$ is a principal ideal domain.
\item An element $r \in R$ is a unit if and only if it satisfies
$\vert r \vert=1$.
\end{itemize}
\end{prop-def}

\begin{definition}
Let $K$ be a field with an absolute value, $R\subset K$ be a Hermite subring, $L$ be a free $R$-module of rank $n$ and 
$b: L \times L \fd K$ be an $R$-symmetric bilinear form. Then
the absolute value of the determinant of the Gram matrix of $b$
is independent of the choice of a base. This real number is the 
discriminant of $b$ and is denoted by
$\left\vert \D(b)\right\vert$.
\end{definition}

In this context, one can generalize Hermite inequality.

\begin{theoreme}[Generalized Hermite inequality]
\label{hermite_generalise}
Let $K$ be a field with an absolute value, $R\subset K$ be a Hermite subring
$L$ be a free $R$-module of rank $n$ and 
$b: L \times L \fd K$ be an $R$-symmetric bilinear form.
Define 
$$ \mu(b) := \min\limits_{x \in L - \{0\}} b(x,x) \pointloin$$ 
Then the following inequality holds:
$$ \mu(b) \leq \Bp{ \frac{1}{1- \rho^2}}^{\frac{n-1}{2}} \left\vert \D(b)\right\vert^{\frac{1}{n}} \qquad,$$
where $\rho$ is the constant attached to $R$.\\
Moreover, if the absolute value of $K$ is non-archimedean, then 
the inequality sharpens to:
$$ \mu(b) \leq \left\vert \D(b)\right\vert^{\frac{1}{n}} \pointloin$$
\end{theoreme}

\begin{exemples} 
\begin{enumerate}
\label{exemple_hermite}
\item When $K=\R$ with its usual absolute value, 
the subring $R=\Z$  is of Hermite for $\rho=\frac{1}{2}$. In this case,
theorem~\ref{hermite_generalise} gives back the 
classical Hermite inequality. 
\item \label{exemple_hermite_fr}Let $q>1$ be a real number and $k$ be a field. We define an absolute value on the field of rational 
functions $K=k(T)$ in the following way.
We set $\vert 0 \vert =0$; and for every non-zero rational function $f$,
we write  $f= \frac{A}{B}$ where $A$ and $B$ are non-zero polynomials in  $k[T]$ and we set
$$ \vert f \vert := q^{\deg A - \deg B} \pointloin$$
This is well defined and one checks that this absolute value 
is a non-archimedean.
The subring $k[T]\subset K$ is of Hermite with constant 
$\rho= \frac{1}{q}$. Indeed, every rational function $f$ is the sum 
of polynomial (its integer part) and of a rational function
whose numerator is of degree strictly less than the degree of the denominator).
\end{enumerate}
\end{exemples}

\subsection{Proof of proposition~\ref{pi0sl_revu}} 

The proof of proposition~\ref{pi0sl_revu} goes by induction on the degree $n$.\\

In the case $n=1$, the matrix $S(T)$ is constant, so there is nothing to prove.\\

Now on, we adopt the conventions of the example~\ref{exemple_hermite} (\ref{exemple_hermite_fr}).
Thus $K=k(T)$ is a field with a non-archimedean norm.
Let $b:k[T]^n \times k[T]^n  \fd k[T]$ be the $K[T]$-symmetric bilinear form
associated to $S(T)$. Since $S(T)$ is non-degenerate, 
we have $\vert \D(b) \vert =1$.
The generalized Hermite inequality implies that
$$ \mu(b) \leq 1 \qquad.$$
Since $S$ has polynomial coefficients,  $\mu(b)$ 
is an integer and thus only two cases are possible:
\begin{itemize}
\item Whether $\mu(b)=1$. This means that there exists a vector $x \in k[T]^n$ such that 
$b(x,x)= \lambda \in (k[T])^{\times}=k^{\times}$.
The Gram matrix of the form  $b$ in 
  $L= \langle x \rangle \oplus \langle x \rangle^{\perp}$,
has the following shape
$$ \begin{bmatrix} 
\lambda & 0 &\dots &0\\
0 & * &\dots& *\\
\vdots & \vdots & \vdots& \vdots\\
 0 & * &\dots& *
\end{bmatrix} \qquad.$$
Applying the induction hypothesis to the restriction of
$b$ to $\langle x \rangle^{\perp}$ is enough to conclude.
 \item Whether $\mu(b)=0$. In this case, there exists a vector $x \in k[T]^n$ such that $b(x,x)=0$. 
Up to renormalisation, we can assume that this $x$
is  indivisible. Since $b$ is non-degenerate, 
there exists a vector $y \in k[T]^n$ such that $b(x,y)=1$. 
The restriction of $b$ to $\langle x,y \rangle$ is then non-degenerate: its Gram matrix in the basis $\{x,y\}$
is of the form
$$ \begin{bmatrix} 
0 & 1 \\
1 & \alpha
\end{bmatrix} \qquad,$$
with  $\alpha \in k[T]$.
We conclude by using the inductive hypothesis on the 
restriction of $b$ to
$\langle x,y \rangle^{\perp}$. 
\end{itemize}
This concludes the proof of proposition~\ref{pi0sl_revu}.
\hfill $\Box$

\section{Additions of rational functions}
\selectlanguage{english}
\label{appendice_comparaison}

The goal of this appendix is to compare the two addition laws
on homotopy classes of endomorphisms of $\pun$: the naive law
denoted $\plus$ defined  in \S\ref{subsec_loi_plus} and the 
$\aun$ law, coming from the fact that $\pun$ is a suspension. 
 More precisely, we are going to prove
 proposition~\ref{prop:plus_oplus_texte}, that is:
\begin{proposition}
\label{prop:plus_oplus}
 The canonical map
$$ \bp{\ppn , \plus} \fd \bp{\pp^\mot, \oplus^\mot}$$
is a monoid morphism.
\end{proposition}

\begin{proof}
Let $g_1$ and $g_2$ be two pointed rational functions. We need to prove that
the image of the rational function $g_1 \plus g_2$ in $\pp^\mot$ is equal to 
$g_1 \oplus^{\mot} g_2$. 
Since the monoid $\ppn$ of naive homotopy classes of rational functions
is generated by its elements of degree $1$ 
(\cf lemma~\ref{lemme_somme_monomes}), 
it is enough to deal with the case when
$g_1$ is of degree $1$. Up to homotopy, one can even assume 
that $g_1$ is of the
form $\frac{X}{a}$ for some $a \in k^\times$.
 For $g_1$ of this form, the formula
for the $\plus$-sum is
$$ \frac{X}{a} \plus g_2= \frac{X}{a}-\frac{1}{a^2} \frac{1}{g_2} \pointloin$$ 

\begin{definition}
 Let $\pun \vtruc \pun$ be the cofiber of the map
$\sph^0=\{0,\infty \}\hookrightarrow \pun \udisj \pun$. 
(Equivalently, $\pun \vtruc \pun$ is the union of two copies of 
$\pun$ with the point 
$0$ in the first copy identified with the point $\infty$ in the second one). The base point is taken at $\infty$.
As $\pun$ is up to homotopy an \emph{unreduced} suspension, there is a canonical map\footnote{
An explicit model for this map is for example the cofiber map
\emph{$\pun \fd \sura{\pun}{\bp{\pun - \{ 0,\infty\}}}$}.
} in the 
$\aun$-homotopy category:
$$ \pun \quad \xrightarrow{\quad \tilde{\nabla} \quad} \quad \pun \vtruc \pun$$ 
$$ { \vcenter{\xymatrix{
 &*=0{\overset{\infty}{{\phantom{x}}}} \ar@{-} '[ddr] [dddd] \ar@{-} '[ddl] [dddd]&  \\
\\
 *=0{} \ar@{-}^-{}[rr]& & *=0{}&\\
\\
& *=0{\underset{0}{\phantom{X}}}& &
}} \xrightarrow{\quad \widetilde{\nabla} \quad}
\qquad
 \vcenter{\xymatrix{
&*=0{\overset{\infty}{{\phantom{x}}}} \ar@{-} '[dr] [dd] \ar@{-} '[dl] [dd]&  \\
 *=0{} \ar@{-}[rr]& & *=0{}&\\
& *=0{\underset{\infty}{\overset{0}{\phantom{X}}}} 
\ar@{-} '[dr] [dd] \ar@{-} '[dl] [dd]  &\\
*=0{} \ar@{-}[rr]& & *=0{}&\\
& *=0{\underset{0}{\phantom{X}}}
}}} \pointloin
 $$
\end{definition}

\begin{lemme}
 Let $g: \pun \fd \pun$ be a pointed rational function and let $f$ be the pointed rational function $f=\frac{X}{a}-\frac{1}{g}$.
One has $f^{-1}(\infty)=\{\infty\}\coprod g^{-1}(0)$ and we denote
$\overline{f}$ the induced map between the cofibers:
$$\overline{f}:  \sura{\pun}{\pun-\bp{\{\infty\}\udisj g^{-1}(0)}} \fd 
\sura{\pun}{\pun-\{ \infty \}} \pointloin$$
 Then the following diagram
of spaces (in the sense of Morel and Voevodsky \cite{MV}) 
$$\vcenter{\xymatrix{ \pun \ar[rd] \ar[rr]^-{f} \ar[ddd]_-{\tilde{\nabla}}  & & \pun \ar[d]^-{\approx} \\
& \sura{\pun}{\pun-\bp{\{\infty\}\udisj g^{-1}(0)}} \ar[r]^-{\overline{f}}  & \sura{\pun}{\pun-\{ \infty \}} \\
 & \Bp{\sura{\pun}{\pun-\{\infty\}}} \vtrucp{*}{*} \Bp{\sura{\pun}{\pun-{g^{-1}(0)}} }  \ar@{=}[u] &\\
\pun \vtruc \pun  \ar[ru] \ar[rr]_-{\frac{X}{a}
 \vtruc -\frac{1}{g}} &
  &  \pun \ar[uu]_-{\approx}
}}$$
homotopy commutes.
\end{lemme}
  
\begin{proof}
 The main point is to prove that the two following diagrams:
\emph{$$ \vcenter{
\xymatrix{
\pun \ar[d] \ar[rr]^-{\frac{X}{a}} & & \pun \ar[d] \\
\sura{\pun}{\pun-\{ \infty \}} \ar@{^(->}[r] &
 \sura{\pun}{\pun-\bp{\{\infty\}\udisj g^{-1}(0)}} \ar[r]^-{\overline{f}} &  \sura{\pun}{\pun-\{ \infty \}}
}}$$} and 
\emph{$$ \vcenter{
\xymatrix{
\pun \ar[d] \ar[rr]^-{g} & & \pun \ar[d] \\
\sura{\pun}{\pun- g^{-1}(0)} \ar@{^(->}[r] &
 \sura{\pun}{\pun-\bp{\{\infty\}\udisj g^{-1}(0)}} \ar[r]^-{\overline{f}} &  \sura{\pun}{\pun-\{ \infty \}}
}}$$}
homotopy commute. The proof is the same in both cases, 
so we give the details only for the last diagram.\\
According to the homotopy purity theorem of Morel and Voevodsky
 (\cf \cite{MV}, theorem~2.23, p. 115),  the space
$\sura{\pun}{\pun-g^{-1}(0)}$ is homotopy equivalent to the Thom space
of the normal bundle of the closed immersion $g^{-1}(0) \inj \pun$.
Since $g^{-1}(0) \cap \{\infty\} = \emptyset$, one has a homotopy 
equivalence 
\emph{$\sura{\aun}{\aun-g^{-1}(0)} \xrightarrow{\approx} \sura{\pun}{\pun-g^{-1}(0)}$} such that the following diagram commutes:
\emph{$$\vcenter{\xymatrix@C=10ex{ \aun \ar[d] \ar[r]^-{f_{\vert \aun}}& \pun \ar[d] \\
\sura{\aun}{\aun-g^{-1}(0)} \ar[r]^-{\overline{f_{\vert \aun}}} \ar[d]_-{\approx}& \sura{\pun}{\pun- \{ \infty \}} \\
\sura{\pun}{\pun-g^{-1}(0)} \ar[ru]_-{\overline{f}}
}} \pointloin$$}
But now, there is a naive homotopy of morphisms of pairs:
\appsnbis{\bp{\aun, \aun-g^{-1}(0)} \times \aun}{\bp{\pun,\pun-\{\infty\}}}{(X,\quad T)}{T\frac X a-\frac{1}{g(X)}}
between $f_{\vert \aun}$ and $(-\frac 1 g)_{\vert \aun}$, which implies  the homotopy commutativity of the diagram considered.\\
 This concludes the proof of the lemma. 
\end{proof}

By the previous lemma, to prove proposition~\ref{prop:plus_oplus} one has to express the composite map
$$\xymatrix{\pun \ar[r]^-{\tilde{\nabla}}&\pun \vtruc \pun   \ar[rr]^-{\frac{X}{a}
 \vtruc -\frac{1}{g}} 
  & & \pun}$$
as a sum in the group $\pp^\mot$. The result 
is given in the next lemma, which is stated in the following context:
\begin{itemize}
\item Let $Y$ be the non-reduced suspension of a pointed space. The base point in $Y$ is taken at $\infty$;
 the other distinguished point is denoted $0$.
\item Let $(Z;z_0)$ be an $\aun$-connected pointed space.
\item Let $g_1,g_2: Y \fd Z$ be two maps such that
$g_1(0)=g_2(\infty)$. We assume that $g_1$ is pointed, that is 
to say that $g_1(\infty)=z_0$.
\item Let $f$ be the pointed map 
$(g_1 \vtruc g_2) \rond \tilde{\nabla}: Y \fd Z$.
\item Let $H: Y \times \aun \fd Z$ be a naive homotopy 
from $g_2$ to a pointed map, say $g_3:Y \fd Z$.
\item Let $\alpha:\aun \fd Z$ be the path $T \mapsto H(\infty, T)$ in $Z$ 
from $g_2(\infty)$ to $z_0$. 
\item Let $\beta: \aun \fd Z$ be the image through $g_1$ of the  canonical path in $Y$ from $\infty$ to $0$;  $\beta$ is thus a path in $Z$ from $z_0$ to $g_1(0)=g_2(\infty)$.
\item Let $\gamma$ be the concatenation of the paths $\beta$ and $\alpha$; $\gamma$ is thus a loop in $Z$ at $z_0$.
\end{itemize}
\begin{lemme}
\label{lemme:somme}
In the group $[Y, Z]^\mot$, one has the  identity:
$$ \gamma \cdot f = (\gamma  \cdot g_1) \oplus^\mot g_3 \pointloin
$$
Above here, the dots ``$\cdot$'' denote the action of an element of $\pi_1^\mot(Z;z_0)$ on $[Y,Z]^\mot$.
\end{lemme}

\begin{proof}
It is a consequence of the following facts:
\begin{itemize}
\item Up to homotopy, $Y$ can be replaced by $\aun \vtruc Y$, pointed at $1\in \aun$. The element $\gamma \cdot f$ is represented by the map 
$$   \vcenter{\xymatrix@C=2ex@R=2ex{
& *=0{\aun \vtruc Y} & &&& \ar@{<-}[llll]+<5ex,0ex>_-{\gamma \vtruc f}Z \\
 &*=0{ \hspace{-2ex} _{\scriptsize 1}} \ar@{-}_-{\mathbf{A}^1}[dd] & &  \quad &\\
& & &  & &\\
&*=0{\overset{\hspace{-2ex} 0}{\underset{ \scriptsize \infty}{{\phantom{_1}}}}} \ar@{-} '[dr] [dd] \ar@{-} '[dl] [dd] \ar[u]& & & &Z  \ar@{<-}@/^/[llld]+<1ex,0ex>^{f}  \ar@{<-}@/_/ [lllu]+<1ex,0ex>_-{\gamma}\\
 *=0{} \ar@{-}[rr]& & *=0{} &   &\\
&  *=0{\underset{0}{\phantom{_1}}} &
}} \pointloin$$
\item Up to homotopy, $Y \vtruc Y$ can be replaced by $ Y \vtrucp{0}{0} \aun \vtrucp{1}{\infty} Y $, the map $Y \vtruc Y \xrightarrow{g_1 \vtruc g_2} Z$ is then homotopic to the map
$$   \vcenter{\xymatrix@C=2ex@R=2ex{
& *=0{Y \vtrucp{0}{0} \aun \vtrucp{1}{\infty} Y} &&&  \qquad \qquad &\ar@{<-}[llll]+<10ex,0ex>_-{g_1 \vtrucp{0}{0} \alpha \vtrucp{1}{\infty} g_3} Z\\
&*=0{{\overset{\infty}{{\phantom{_1}}}}} \ar@{-} '[dr] [dd] \ar@{-} '[dl] [dd]& & \\
 *=0{} \ar@{-}[rr]& & *=0{}&   &\\
&  *=0{\overset{0}{\underset{\hspace{-2ex}\raisebox{0,5ex}{\tiny 0}}{\phantom{X}}}} \ar@{-}^-{\mathbf{A}^1}[dd] \ar[d]& \\
& & & &   & Z \ar@{<-}[ll]+<0ex,0ex>_-{\alpha} \ar@{<-}@/_1pc/[lluu]+<0ex,0ex>_{g_1} \ar@{<-}@/^1pc/[lldd]+<0ex,0ex>^{g_3}\\
 &*=0{\overset{\hspace{-2ex} 1}{\underset{\infty}{{\phantom{_1}}}}} \ar@{-} '[dr] [dd] \ar@{-} '[dl] [dd]  &\\
*=0{} \ar@{-}[rr]& & *=0{}&  &\\
& *=0{\underset{0}{\phantom{X}}}
}}$$
\item The restriction of the map 
$$\aun \vtruc Y \vtrucp{0}{0} \aun \vtrucp{1}{\infty} Y \xrightarrow{ \gamma \vtruc f \vtrucp{0}{0} \alpha \vtrucp{1}{\infty} g_3} Z$$
 to
$\aun \vtrucp{0}{0} \aun \vtrucp{1}{0} \aun$ (the $\aun$ in the middle is the domain of $\beta$) is the concatenation of $\gamma^{-1}$ and $\gamma$ and is thus null-homotopic.
 So, up to homotopy, the map  
$\gamma \vtruc f \vtrucp{0}{0} \alpha \vtrucp{1}{\infty} g_3$ 
factors through the   cofiber
$\sura{\bp{\aun \vtruc Y \vtrucp{0}{0} \aun \vtrucp{1}{\infty} Y}}{\bp{\aun \vtrucp{0}{0} \aun \vtrucp{1}{0} \aun
}} \approx Y\vee Y$.
\item The following composite map  
\emph{$$ \vcenter{\xymatrix{ Y  \ar@{->}[r]^-{\approx} \ar@{-->}[d]& \aun \vtruc Y \ar[r]^-{\id\vtruc \tilde{\nabla}} &
{\aun \vtruc Y \vtruc Y} \ar[d]^-{\approx} \\
 Y\vee Y & \sura{\bp{\aun \vtruc Y \vtrucp{0}{0} \aun \vtrucp{1}{\infty} Y}}{\bp{\aun \vtrucp{0}{0} \aun \vtrucp{1}{0} \aun }} \ar[l]_-{\approx}& {\aun \vtruc Y \vtrucp{0}{0} \aun \vtrucp{1}{\infty} Y} \ar[l] &   }}
$$}
is  a model for the co-diagonal $Y \xrightarrow{\nabla} Y \vee Y$.
\end{itemize}
\end{proof}

In our situation, there is a ``universal'' homotopy between $-\frac{1}g$ and $g$, which is given by composition at the target of $g$ with a naive homotopy between
between $- \frac 1 X$ and $X$.
One can choose  
$$ 
{\frac{(-T^2+ 2T)X -(T^3-3T^2+T+1)}{(-T+1)X + (-T^2+2T)}}
$$ 
as an example of such a homotopy.

\begin{lemme}
\label{lemme:action}
 Let $\alpha: \aun \fd \pun$ be the path $T\mapsto [T : 1-T^2]$,
 $\beta: \aun \fd \pun$ be the path $ T \mapsto [1-T : aT]$ and $\gamma$ 
be the loop given by concatenation of $\alpha$ and $\beta$. 
Then for any
pointed rational function $F:\pun \fd \pun$, 
one has the identity
$$\gamma \cdot F= a^2 F$$
 in  $\pp^\mot$.
\end{lemme}

\begin{proof}
The projective line $\pun$ is homotopy equivalent to the homogeneous space
\emph{${\sura{\SL_2}{\gm}}$}. So any path to $\pun$ lifts to a path in $\SL_2$. 
Moreover, up to homotopy, one can use the product in $\SL_2$ to compose paths.
  More precisely, let $\tilde{\beta}: \aun \fd \SL_2$ be a lift of 
$\beta$ such that $\tilde{\beta}(0)=\id$ and let $\tilde{\alpha}: \aun \fd \SL_2$ be a lift of $\alpha$ such that $\tilde{\alpha}(0)=\tilde{\beta}(1)$. 
Then a lift $\tilde{\gamma}$ of the loop $\gamma$ is 
\app{\tilde{\gamma}}{\aun}{\SL_2}{T} 
 {\tilde{\beta}(T) \cdot \tilde{\beta}(1)^{-1} \cdot \tilde{\alpha}(T)} 
Let $c_{\infty}:\aun \fd \pun$ be the constant path at $\infty \in \pun$.
The action of $\SL_2$ on $\pun$ can be used to define a 
naive homotopy
between $\gamma \vtruc F$ and $c_{\infty} \vtruc (\tilde{\alpha}(1) \cdot F)$, namely\footnote{
Below, $H\in \aun$ is the parameter of the homotopy
and $T\in \aun$ belongs to the domain of $\gamma$.}
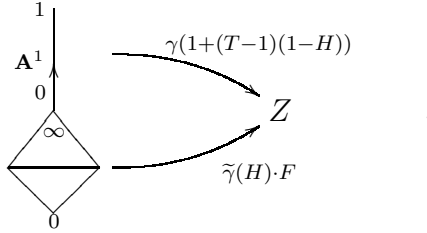
$$\vcenter{\xymatrix@C=2ex@R=2ex{
 &*=0{ \hspace{-2ex} _{\scriptsize 1}} \ar@{-}_-{\mathbf{A}^1}[dd] & &  \quad &\\
& & &  & &\\
&*=0{\overset{\hspace{-2ex} 0}{\underset{ \scriptsize \infty}{{\phantom{_1}}}}} \ar@{-} '[dr] [dd] \ar@{-} '[dl] [dd] \ar[u]& & & &Z  \ar@{<-}@/^/[llld]+<1ex,0ex>^(.5){\qquad {\tilde{\gamma}(H)\cdot F}}  \ar@{<-}@/_/ [lllu]+<1ex,0ex>_-{\qquad \gamma(1+(T-1)(1-H))}\\
 *=0{} \ar@{-}[rr]& & *=0{} &   &\\
&  *=0{\underset{0}{\phantom{_1}}} &
}} \pointloin $$
Explicitly, one can take 
  $\tilde{\alpha}(T)=\begin{bmatrix}
a (-T^2  + 2  T)& (-2 +\frac{1}{a})  T^2+ (4 -  \frac{1}{a})  T - \frac{1}{a}\\ a(- T +1)& (-2  + \frac{1}{a})  T + 2
\end{bmatrix}$ and $\tilde{\beta}(T)= 
\begin{bmatrix}
 1-T & -\frac T a \\ aT & 1+T
\end{bmatrix}
$.
So $\tilde{\gamma}(1)=\tilde{\alpha}(1)=\begin{bmatrix}  
a & 2 - \frac 1a \\ 0 & \frac 1 a
\end{bmatrix}
$. The result follows since $\tilde{\gamma}(1) \cdot F=a^2F +
2-\frac 1a$ is 
canonically homotopic to
$a^2 F$.  
\end{proof}

Lemma~\ref{lemme:somme} and lemma~\ref{lemme:action} together imply  proposition~\ref{prop:plus_oplus}. Indeed, for every $a \in k^\times$ and for every pointed rational function $g$
one has the identity in $\pp^\mot$
$$ aX \plus g= aX- \frac{a^2}{g}=a^2\Bp{\frac{X}{a}-\frac{1}{g}} = aX \oplus^\mot g\pointloin$$
\end{proof}

\begin{remarque}
 It is likely that the preceding method could be refined to 
show that for every positive integer $n$ the diagram of spaces
$$\vcenter{\xymatrix{
\F_1 \times \F_n \ar[r]^-{\plus} \ar@{^(->}[d] & \F_{n+1} \ar@{^(->}[d] \\
\Omega^{\pun}_1 \pun \times \Omega^{\pun}_n \pun \ar[r]^-{\oplus^\mot} & \Omega^{\pun}_{n+1} \pun
}} $$ 
homotopy commutes.
\end{remarque}

\section{Hankel matrices}
\label{sec_hankel}

Let $n$ be a positive integer. We described in \S\ref{subsec_bezout}
a way to produce a non-degenerate \nn{} symmetric matrix
out of a pointed rational function. However, for dimension reasons,
it appears that for $n\geq3$, not any non-degenerate \nn{} symmetric matrix
is the Bézout form of some rational function.
Indeed, one observes (\cf lemma~\ref{inverse_bezout_hankel})
 that the inverse of the Bézout form of 
rational function is \emph{Hankel}, that is has constant value 
along anti-diagonals, \cf definition~\ref{def_hankel}.
It turns out that this necessary condition is also sufficient.
The goal of this appendix is to review this old story with our
conventions.

\begin{definition}
\label{def_hankel}
A symmetric matrix $S$ is of Hankel if the value of its 
entries $s_{p,q}$  depends only on $p+q$. For every integer 
$n$, let $\mc{H}_n$ be the scheme of non-degenerate
 \nn{} Hankel matrices;  $\mc{H}_n$ is thus a closed subscheme of  $\mc{S}_n$ of dimension $2n-1$.
\end{definition}

\begin{lemme}
\label{inverse_bezout_hankel}
Let $R$ be a ring, $n$ be a positive integer and $f=\frac{A}{B} \in \F_n(R)$ be a rational function.
 Then the matrix  $\bez_n(f)^{-1}$ is Hankel.
\end{lemme}

\begin{proof}
Consider the universal case when the coefficients of
$A$ and $B$ are indeterminates, and let $K$ be the algebraic closure of the quotient field of  $R=\Z[a_i,b_j]$.\\
The polynomial $A$ is split in $K$. Let 
$(\alpha_i)_{1\leq i \leq n}$ be the roots of $A$ in $K$ and let $E$ be the étale $K$-algebra $\sur{K[X]}{(A)}$.
We are going to show that 
$\bez_n(f)^{-1}$ 
is the matrix in the canonical basis
$(1,X, \dots, X^{n-1})$ of $E$,
of the following trace form
\appsn{E \times E}{K}{(P,\, Q)}{\mathrm{tr}_{E/K}\left(\frac{PQ}{A^\prime B}\right)}
The generic entry of this matrix is 
 $\mathrm{tr}_{E/K}\bp{\frac{X^{p+q}}{A^\prime B}}$
which proves that $\bez_n(f)^{-1}$  is indeed Hankel.\\

Let $V$ be the Vandermonde matrix associated to the roots of $A$, that is
$$ V:= \begin{bmatrix} 1 &1 & \dots & 1\\
\alpha_1 & \alpha_2 & \dots & \alpha_n \\
\alpha_1^2& \alpha_2^2 & \dots & \alpha_n^2 \\
\dots &      \dots &    \dots &  \dots \\
\alpha_1^{n-1}& \alpha_2^{n-1} & \dots & \alpha_n^{n-1} 
\end{bmatrix} \qquad.$$
By construction of  $\bez_n(A,B)$, we have the following identity
 $$ \transp V \cdot \bez_n(A,B)\cdot V = \begin{bmatrix}
\delta_{A,B}(\alpha_p,\alpha_q) 
\end{bmatrix} \qquad,$$
where $\delta_{A,B}$ was introduced in the definition of $\bez_n(A,B)$.
For every integers $1 \leq p,q \leq n$, one has 
$$ \delta_{A,B}(\alpha_p,\alpha_q) = \begin{cases} A^\prime(\alpha_p) B(\alpha_p) & \text{ whenever } p=q\\
0 & \text{ otherwise}
\end{cases} \qquad.$$
This shows that $\bez_n(f)^{-1}$ is indeed the matrix of the
previous trace form.
\end{proof}

\begin{remarque}
\label{remarque_developpement}
We keep the notations introduced in the proof of lemma~\ref{inverse_bezout_hankel}.
Let also $\rho$ be the linear form \emph{residue}
\app{\rho}{E}{K}{X^i}{\begin{cases} 
1 & \text{ if } i=n-1\\
0 & \text{ otherwise}
\end{cases}}
By a classical identity due to Euler (see \cite{serre_corps_locaux}, chapitre~2, $\mathrm{n}^{\rond}$6, lemme~2 for example), we have 
$$ \forall P \in E, \qquad \rho(P)=\tr_{E/K}\left(\frac{P}{A^\prime}\right) \qquad.$$
The matrix $\bez_n(A,B)^{-1}$ is thus the matrix of the 
symmetric bilinear form
$(P,Q) \mapsto \rho\left(\frac{PQ}{B}\right)= \rho(PQV)$, 
where $V$ is a polynomial associated to a 
Bézout relation $AU+BV=1$.\\
One thus deduces the following algorithm to compute
the Hankel matrix
$\bez_n(A,B)^{-1}$.
One formally expands in $R[X,X^{-1}]]$ the rational function $\frac{V}{A}$:
\begin{equation}
\label{developpement}
 \frac{V}{A}=s_1 X^{-1} + s_2 X^{-2} + \dots + s_{2n-1}X^{-(2n-1)} + s_{2n} X^{-2n} +\mathrm{O}(X^{-(2n+1)}) \qquad.
\end{equation}
The matrix  $\bez_n(A,B)^{-1}$ is then given by
$$\bez_n(A,B)^{-1}= \begin{bmatrix} 
s_1 & s_2 & s_3 &\dots & s_n \\
s_2 & s_3 & \dots &\dots & s_{n+1}\\
s_3 &\dots & \dots & \dots & \dots \\
\dots & \dots & \dots & \dots & \dots \\
s_n & s_{n+1} & \dots & \dots &s_{2n-1}
\end{bmatrix} \qquad.$$
\end{remarque}

\begin{definition} Let $R$ be a ring and let $n$ be a positive integer.
\begin{itemize}
\item For all $(A,B) \in \F_n(R)$, there exists a unique pair of polynomials $(U_1, V_1)$ with  $\deg U_1=n-1$, $\deg V_1\leq n-1$
and such that 
$$ AU_1 + BV_1 = X^{2n-1} \qquad.$$ 
One defines a scheme morphism
$\phi_{n}: \F_n \fd \aun$ 
by associating to the pair $(A,B)$ the opposite of the coefficient of
 $X^{n-1}$ in $V_1$. 
(Note that the scalar $\phi_{n}(A,B)$ just defined is also 
equal to the opposite of the coefficient  $s_{2n}$ of $X^{-2n}$
in the formal expansion (\ref{developpement}) of  $\frac{V}{A}$).
\item Let $\hank_n$ be the morphism of schemes
\app{\hank_n}{\F_n}{\mc{H}_n}{(A,B)}{\bez_n(A,B)^{-1}}
\end{itemize}
\end{definition}

\begin{proposition}
\label{proposition_isom_hankel}
For every integer $n$, the morphism
$$ \F_n \xrightarrow{\hank_n \times\, \phi_{n}}  \mc{H}_n \times \aun$$
is a $\ga$-equivariant isomorphism of schemes.
\end{proposition}

\begin{proof}
The morphism was checked to be equivariant in \S\ref{n2}.
Moreover, remark~\ref{remarque_developpement} 
suggests a way to produce the inverse morphism.
Let
\app{\psi_n}{\mc{H}_n \times \aun }{\F_n}{(H,v)}{(A,B)}
where $A$ is the degree $n$ monic polynomial whose coefficients $a_i$ 
are given by 
$$\begin{bmatrix}
a_0 \\ a_1 \\ \dots \\ a_{n-2}\\a_{n-1}
\end{bmatrix} = H^{-1} \begin{bmatrix}
-s_{n+1} \\ -s_{n+2} \\ \dots \\ -s_{2n-1} \\ v
\end{bmatrix} \qquad,$$
and where $B$ is the unique polynomial of degree $\leq n-1$
given by the Bézout  relation  $AU+BV=1$ for
$V=\mathrm{E}((\sum_{i=1}^{n} s_i X^{-i})A)$, $\mathrm{E}$ denoting the
integer part.\\
By construction, $\psi_n$ is inverse to $\hank_n \times \phi_{n}$.
\end{proof}

\begin{remarque}
 The scheme $\F_n$ is a $\ga$-torsor over the base 
$\mc{H}_n$, which is affine. Thus it is not surprising that this
torsor is trivial. 
\end{remarque}

\begin{remarque}
 Proposition~\ref{proposition_isom_hankel} allows to reformulate our
 results in terms of Hankel matrices. For example, there is a canonical 
monoid structure denoted $\plus$ on the set $\udisj\limits_{n\geq 0}(\pN\mc{H}_n)(k)$
such that the map 
\appsn{\mc{H}_n}{\mc{S}_n}{H}{H^{-1}}
induces a monoid isomorphism
$$ \Bp{\udisj\limits_{n \geq0}(\pN\mc{H}_n)(k), \plus} \xrightarrow{\isom} \Bp{\udisj\limits_{n \geq0}(\pN\mc{S}_n)(k), \oplus} \pointloin$$
\end{remarque}

\def\cprime{$'$}
\providecommand{\bysame}{\leavevmode ---\ }
\providecommand{\og}{``}
\providecommand{\fg}{''}
\providecommand{\smfandname}{\&}
\providecommand{\smfedsname}{\'eds.}
\providecommand{\smfedname}{\'ed.}
\providecommand{\smfmastersthesisname}{M\'emoire}
\providecommand{\smfphdthesisname}{Th\`ese}

\end{document}